\documentstyle[12pt,leqno]{article}
\def\R{{\bf R}}
\def\C{{\bf C}}

\def\P{{\bf P}}

\begin{document}
\baselineskip=0.52cm

 $\:$  $\:$ $\:$  $\:$ $\:$  $\:$ $\:$  $\:$ $\:$  $\:$ 
 $\:$  $\:$ $\:$  $\:$ $\:$  $\:$ $\:$  $\:$ $\:$  $\:$ 
 $\:$  $\:$ $\:$  $\:$ $\:$  $\:$ $\:$  $\:$ $\:$  $\:$ 
 $\:$  $\:$ $\:$  $\:$ $\:$  $\:$ $\:$  $\:$ $\:$  $\:$ 
 $\:$  $\:$ $\:$  $\:$ $\:$  $\:$ $\:$  $\:$ $\:$  $\:$ 
 $\:$  $\:$ $\:$  $\:$ $\:$  $\:$ $\:$  $\:$ $\:$  $\:$ 
 $\:$  $\:$ $\:$  $\:$ $\:$  $\:$ $\:$  $\:$ $\:$  $\:$ 
 $\:$  $\:$ $\:$  $\:$ $\:$  $\:$ $\:$  $\:$ $\:$  $\:$ 
 $\:$  $\:$ $\:$  $\:$ $\:$  $\:$ $\:$  $\:$ $\:$  $\:$ 
 $\:$  $\:$ $\:$  $\:$ $\:$  $\:$ $\:$  $\:$ $\:$  $\:$ 
 $\:$  $\:$ $\:$  $\:$ $\:$  $\:$ $\:$  $\:$ $\:$  $\:$ 
 $\:$  $\:$ $\:$  $\:$ $\:$  $\:$ $\:$  $\:$ $\:$  $\:$ 
 $\:$  $\:$ $\:$  $\:$ $\:$  $\:$ $\:$  $\:$ $\:$  $\:$ 
 $\:$  $\:$ $\:$  $\:$ $\:$  $\:$ $\:$  $\:$ $\:$  $\:$ 
 $\:$  $\:$ $\:$  $\:$ $\:$  $\:$ $\:$  $\:$ $\:$  $\:$ 
 $\:$  $\:$ $\:$  $\:$ $\:$  $\:$ $\:$  $\:$ $\:$  $\:$ 
 $\:$  $\:$ $\:$  $\:$ $\:$  $\:$ $\:$  $\:$ $\:$  $\:$ 
 $\:$  $\:$ $\:$  $\:$ $\:$  $\:$ $\:$  $\:$ $\:$  $\:$ 
 $\:$  $\:$ $\:$  $\:$ $\:$  $\:$ $\:$  $\:$ $\:$  $\:$ 
 $\:$  $\:$ $\:$  $\:$ $\:$  $\:$ $\:$  $\:$ $\:$  $\:$ 
 $\:$  $\:$ $\:$  $\:$ $\:$  $\:$ $\:$  $\:$ $\:$  $\:$ 
 $\:$  $\:$ $\:$  $\:$ $\:$  $\:$ $\:$  $\:$ $\:$  $\:$ 
 $\:$  $\:$ $\:$  $\:$ $\:$  $\:$ $\:$  $\:$ $\:$  $\:$ 
 $\:$  $\:$ $\:$  $\:$ $\:$  $\:$ $\:$  $\:$ $\:$  $\:$ 
 $\:$  $\:$ $\:$  $\:$ $\:$  $\:$ $\:$  $\:$ $\:$  $\:$ 
 $\:$  $\:$ $\:$  $\:$ $\:$  $\:$ $\:$  $\:$ $\:$  $\:$ 
  
 $\:$  $\:$ $\:$  $\:$ $\:$  $\:$ $\:$  $\:$ $\:$  $\:$ 
 $\:$  $\:$ $\:$  $\:$ $\:$  $\:$ $\:$  $\:$ $\:$  $\:$ 
 $\:$  $\:$ $\:$  $\:$ $\:$  $\:$ $\:$  $\:$ $\:$  $\:$ 
 $\:$  $\:$ $\:$  $\:$ $\:$  $\:$ $\:$  $\:$ $\:$  $\:$ 
 $\:$  $\:$ $\:$  $\:$ $\:$  $\:$ $\:$  $\:$ $\:$  $\:$ 
 $\:$  $\:$ $\:$  $\:$ $\:$  $\:$ $\:$  $\:$ $\:$  $\:$ 
 $\:$  $\:$ $\:$  $\:$ $\:$  $\:$ $\:$  $\:$ $\:$  $\:$ 
 $\:$  $\:$ $\:$  $\:$ $\:$  $\:$ $\:$  $\:$ $\:$  $\:$ 
 $\:$  $\:$ $\:$  $\:$ $\:$  $\:$ $\:$  $\:$ $\:$  $\:$ 
 $\:$  $\:$ $\:$  $\:$ $\:$  $\:$ $\:$  $\:$ $\:$  $\:$ 
 $\:$  $\:$ $\:$  $\:$ $\:$  $\:$ $\:$  $\:$ $\:$  $\:$ 
 $\:$  $\:$ $\:$  $\:$ $\:$  $\:$ $\:$  $\:$ $\:$  $\:$ 
 $\:$  $\:$ $\:$  $\:$ $\:$  $\:$ $\:$  $\:$ $\:$  $\:$ 
 $\:$  $\:$ $\:$  $\:$ $\:$  $\:$ $\:$  $\:$ $\:$  $\:$ 
 $\:$  $\:$ $\:$  $\:$ $\:$  $\:$ $\:$  $\:$ $\:$  $\:$ 
 $\:$  $\:$ $\:$  $\:$ $\:$  $\:$ $\:$  $\:$ $\:$  $\:$ 
 $\:$  $\:$ $\:$  $\:$ $\:$  $\:$ $\:$  $\:$ $\:$  $\:$ 
 $\:$  $\:$ $\:$  $\:$ $\:$  $\:$ $\:$  $\:$ $\:$  $\:$ 
 $\:$  $\:$ $\:$  $\:$ $\:$  $\:$ $\:$  $\:$ $\:$  $\:$ 
 $\:$  $\:$ $\:$  $\:$ $\:$  $\:$ $\:$  $\:$ $\:$  $\:$ 
 $\:$  $\:$ $\:$  $\:$ $\:$  $\:$ $\:$  $\:$ $\:$  $\:$ 
 $\:$  $\:$ $\:$  $\:$ $\:$  $\:$ $\:$  $\:$ $\:$  $\:$ 
 $\:$  $\:$ $\:$  $\:$ $\:$  $\:$ $\:$  $\:$ $\:$  $\:$ 
 $\:$  $\:$ $\:$  $\:$ $\:$  $\:$ $\:$  $\:$ $\:$  $\:$ 
 $\:$  $\:$ $\:$  $\:$ $\:$  $\:$ $\:$  $\:$ $\:$  $\:$ 
 $\:$  $\:$ $\:$  $\:$ $\:$  $\:$ $\:$  $\:$ $\:$  $\:$

\begin{center}{\large\bf
ON THE LOCAL MEROMORPHIC EXTENSION  \\
OF CR MEROMORPHIC MAPPINGS
\vspace{0.16cm}}
\end{center}

\begin{center}
{\bf Jo\"el Merker} and
{\bf Egmont Porten}
\end{center}

\bigskip
\bigskip
\begin{center}
{\large Laboratoire d'Analyse, Topologie et Probabilit\'es, UMR 6632, \\
CMI, 39 rue Joliot Curie, 13453 Marseille Cedex 13, France \\
E-Mail address: merker@gyptis.univ-mrs.fr}

\bigskip
{\large Max-Planck-Gesellschaft, Humboldt-Universit\"at zu Berlin \\
J\"agerstrasse, 10-11, D-10117 Berlin, Germany \\
E-Mail address: egmont@mathematik.hu-berlin.de}

\bigskip
\bigskip

\begin{minipage}[t]{11cm}{\small
{\sc Abstract.}\footnotemark \  
Let $M$ be a generic CR submanifold in $\C^{m+n}$, $m=\hbox{CRdim} \ M \geq 1$, 
$n=\hbox{codim} \ M \geq 1$, $d=\hbox{dim} \ M = 2m+n$. A CR meromorphic
mapping (in the sense of Harvey-Lawson) is a triple 
$(f,{\cal D}_f, [\Gamma_f])$, where: 1. $f: {\cal D}_f \to
Y$ is a ${\cal C}^1$-smooth mapping defined over a dense open subset
${\cal D}_f$ of $M$ with values in a projective manifold
$Y$; 2. The closure $\Gamma_f$ of its graph in $\C^{m+n} \times Y$ defines
a oriented scarred ${\cal C}^1$-smooth CR manifold of CR dimension $m$ ({\em i.e.} CR
outside a closed thin set) and 3. Such that $d[\Gamma_f]=0$ in the sense of currents. We prove
in this paper that $(f,{\cal D}_f, [\Gamma_f])$ extends meromorphically
to a wedge attached to $M$ if $M$ is everywhere minimal and 
${\cal C}^{\omega}$ (real analytic) or if
$M$ is a ${\cal C}^{2,\alpha}$ globally minimal hypersurface.}
\end{minipage}
\end{center}

\bigskip
\smallskip
\smallskip


\newpage

                         $\:$
                        \bigskip
                         \begin{center}
                         ON THE LOCAL MEROMORPHIC EXTENSION \\
                         OF CR MEROMORPHIC MAPPINGS
                         \end{center}

                         \bigskip
                         \begin{center}
                         {\bf J. Merker} and {\bf E. Porten}
                         \end{center}

\bigskip
\bigskip

Since the works of Tr\'epreau, Tumanov and J\"oricke, extendability properties
of CR functions on a smooth CR manifold $M$ became fairly well understood. In
a natural way, 1. $M$ is seen to be a disjoint union of CR bricks, 
called {\it CR orbits}, each
of which being an immersed CR submanifold of $M$ with the same CR dimension 
as $M$  (\cite{TR}) ; 2. A continuous CR function $f$ on $M$ is CR if and only if
its restriction $f|_{{\cal O}_{CR}}$ is CR on each CR orbit ${\cal O}_{CR}$ (\cite{JO2}, \cite{PO1}, \cite{MP}) ;
3; For each CR orbit ${\cal O}_{CR}$, there exists an analytic wedge
${\cal W}^{an}$ attached to ${\cal O}_{CR}$, {\em i.e.} a conic complex
manifold with edge ${\cal O}_{CR}$ and with 
$\hbox{dim}_{\C} {\cal W}^{an} =\hbox{dim}_{\R} ({\cal O}_{CR})
-\hbox{CRdim} \ M$, such that each continuous CR function
on ${\cal O}_{CR}$ admits a holomorphic extension to ${\cal W}^{an}$
(\cite{TR}, \cite{JO2}, \cite{TU}, \cite{ME1}). 
The technique of FBI tranforms (\cite{TR}) or deformations of
analytic discs (\cite{TU}, \cite{JO1}, \cite{ME1}) brings up the construction 
of the analytic wedges in a semi-local way.

This paper is devoted to the question of meromorphic extension to 
wedges of CR meromorphic functions in the sense of 
Harvey and Lawson (\cite{HL}, see also \cite{SA}).

The classical theorem of Hartogs-Levi states that, if a meromorphic
function is given on a neighborhood ${\cal V}(b\Omega)$ of a bounded
domain $\Omega \subset \subset \C^{m+1}$, $m+1 \geq 2$, then 
it extends meromorphically inside $\Omega$. Using the solution of the
complex Plateau problem, {\em i.e.} attaching holomorphic chains
to maximally complex cycles in the complex euclidean space, Harvey and 
Lawson proved the following Hartogs-Bochner theorem for 
meromorphic maps: If $m+1\geq 3$, any CR mapping $b\Omega \to Y$, with values
in a projective manifold $Y$, extends meromorphically to $\Omega$. The method allows indeterminacies :
a CR meromorphic mapping is defined by Harvey and Lawson as a triple
$(f,{\cal D}_f,[\Gamma_f])$, where $f: {\cal D}_f \to Y$ is a ${\cal C}^1$-smooth
 mapping defined over a dense open subset ${\cal D}_f \subset M=b\Omega$ with 
values in a projective manifold $Y$; the closure $\Gamma_f$ of its graph
in $\C^{m+1} \times Y$ defines a scarred ${\cal C}^1$-smooth CR manifold
of CR dimension $m$ ({\em i.e.} CR outside a closed thin set) and such that
$d[\Gamma_f]=0$ in the sense of currents. The case $m=1$ was open until
Dolbeault and Henkin gave a positive answer for $C^2$ CR mappings $f$ using 
their solution of the boundary problem in $\C\P^n$, $n\geq 2$ 
(\cite{DH}). For continuous $f$ with values in a compact K\"ahler manifold
the second-named author devised
a different proof, relying on the fact that $b\Omega$ is 
a single CR orbit and that the envelope of holomorphy 
of ${\cal V}(b\Omega)$ contains
$\Omega$ (\cite{PO2}, see Section 4).

Recently, Sarkis obtained the analog of the Hartogs-Bochner
theorem for meromorphic maps, allowing indeterminacies (\cite{SA}, see Section 4). The
main idea is to see that the set $\Sigma_f$ of indeterminacies
of $(f,{\cal D}_f,[\Gamma_f])$ is a closed subset with empty
interior of some ${\cal C}^1$-scarred submanifold $\Lambda \subset M=b\Omega$, 
with $\hbox{codim}_M \Lambda =2$, and that $f$ defines an order zero 
CR distribution on $M\backslash \Sigma_f$.
Then the question of CR meromorphic extension is reduced to the local
removable singularities theorems in the 
spirit of J\"oricke (\cite{JO88}, \cite{JO1}, \cite{JO2}). 
We would like to mention
that these removability results were originally impulsed 
by J\"oricke in \cite{JO88} and in \cite{JO2}.

The goal of this article is to push forward meromorphic extension on CR 
manifolds of arbitrary codimension, the analogs of domains being {\it wedges} over
CR manifolds. It seems natural to use the theory of Tr\'epreau-Tumanov in this 
context. Knowing thinness of $\Sigma_f$ (Sarkis) and using wedge removable
singularities theorems (\cite{ME2}, \cite{MP}, \cite{PO1}), we prove in this paper
that a CR meromorphic mapping $(f,{\cal D}_f,[\Gamma]_f)$ extends 
meromorphically to a wedge attached to $M$ if the CR generic manifold $M$ is everywhere minimal
in the sense of Tumanov and real analytic, ${\cal C}^{\omega}$-smooth.
We prove also that
such CR meromorphic mappings extend meromorphically
to a wedge if $M$ is a ${\cal C}^{2,\alpha}$-smooth ($0<\alpha<1$)
hypersurface in $\C^{m+1}$ that is only globally minimal and
we prove the meromorphic extension in 
any codimension if $M$ is everywhere minimal and if the scar set $\hbox{Sc}(\Sigma_f)$ (in fact ${\rm Sc}(\Lambda)$) of the indeterminacy 
set $\Sigma_f$ is of $(d-3)$-dimensional 
Hausdorff measure equal to zero, $d=2m+n=\hbox{dim} \ M$.
These results are parallel to the meromorphic extension theorem obtained
by Dinh and Sarkis for manifolds $M$ with 
nondegenerate vector-valued Levi-form (\cite{DS}).

We refer the reader to Section 4 which plays 
the role of a detailed introduction.

\smallskip
{\it Acknoweledgement.} We are grateful to Professor Henkin who raised
the question. We also wish to address special thanks to Frederic
Sarkis. He has communicated to us the reduction of meromorphic
extension of CR meromorphic mappings to a removable singularity
property and we had several interesting conversations with him.

\bigskip
\noindent
{\bf 1. Currents and scarred manif
olds.}  In this section, we follow
Harvey and Lawson for a preliminary exposition of currents in the CR
category. This material is known, and is recalled here for clarity.
Let ${\cal U}\subset \C^{m+n}$ be an open set.  We shall denote by
${\cal D}^k({\cal U})$ the space of all complex-valued ${\cal
C}^{\infty}$ exterior $k$-forms on ${\cal U}$ with the usual
topology. The dual space to ${\cal D}^k({\cal U})$ will be denoted by
${\cal D}_k'({\cal U})$. We adopt the dual notation ${\cal D}_k'({\cal
U})={\cal D}'{\:\!}^{2(m+n)-k}({\cal U})$ and say that elements of
this space are currents of dimension $k$ and degree $2(m+n)-k$ on
${\cal U}$. In fact, every $k$-dimensional current can be naturally
represented as an exterior $(2(m+n)-k)$-form on ${\cal U}$ with
coefficients in ${\cal D}_{2m+2n}'({\cal U})$.

We let $d: {\cal D}^k({\cal U}) \to {\cal D}^{k+1}({\cal U})$ denote
the exterior differentiation operator and also denote by $d: {\cal
D}_{k+1}'({\cal U}) \to {\cal D}_k'({\cal U})$ the adjoint map ({\em
i.e.} $d: {\cal D}'{\:\!}^{2(m+n)-k-1}({\cal U}) \to {\cal
D}'{\:\!}^{2(m+n)-k}({\cal U})$).

In the following, ${\cal H}^q$, $q\in \R$, $0\leq q \leq 2m+2n$, will
denote Hausdorff $q$-dimensional measure on $\C^{m+n}$. The notation
${\cal H}_{loc}^q(E) < \infty$ for a set $E \subset \C^{m+n}$ means
that, for all compact subsets $K \subset \subset E$, ${\cal H}^q(K) <
\infty$.

We have the Dolbeault decomposition ${\cal D}^k({\cal
U})=\oplus_{r+s=k} {\cal D}^{r,s}({\cal U})$ and its dual
decomposition ${\cal D}_k'({\cal U}) = \oplus_{r+s=k} {\cal
D}_{r,s}'({\cal U})$ (or ${\cal D}'{\:\!}^{2(m+n)-k}({\cal
U})=\oplus_{r+s=k} {\cal D}'{\:\!}^{m+n-r,m+n-s}({\cal U}))$. A
current in ${\cal D}_{r,s}'({\cal U})= {\cal
D}'{\:\!}^{m+n-r,m+n-s}({\cal U})$ is said to have bidimension $(r,s)$
and bidegree $(m+n-r,m+n-s)$. Given a current $T\in {\cal D}_k'({\cal
U})$, we will denote the components of $T$ in the space ${\cal
D}_{r,s}'({\cal U})={\cal D}'{\:\!}^{m+n-r,m+n-s}({\cal U})$ by
$T_{r,s}$ of $T^{m+n-r,m+n-s}$ : the subscripts refer to bidimension
and the superscripts to bidegree. Thus
$${\cal D}_k'({\cal U}) \ni T =\sum_{r+s=k} T_{r,s} = \sum_{r+s=k}
T^{m+n-r,m+n-s}.$$

Let $M$ be an oriented $d$-dimensional manifold of
class ${\cal C}^1$ in ${\cal U}$ with ${\cal H}_{loc}^d(M) <
\infty$ (we refer the reader to the paragraph
before Proposition 5.7 for a presentation of Hausdorff measures). 
Then $M$ defines a current $[M] \in {\cal D}_d'({\cal U})$,
called the current of integration on $M$, by $[M] (\varphi) =\int_M
\varphi$, for all $\varphi \in {\cal D}^q({\cal U})$. Furthermore,
$d[M]=0$ if $b M =\emptyset$ by Stokes' formula, in particular if $M$
is a closed submanifold of ${\cal U}$. An obvious remark is that
$[M]=[M\backslash \sigma]$, for all closed sets $\sigma \subset {\cal
U}$ with ${\cal H}^d(\sigma)=0$. For example, pure $d$-dimensional
real or complex analytic sets $\Psi \subset {\cal U}$ have a geometric
decomposition in a regular and a singular part, $\Psi=
\hbox{Reg}(\Psi) \cup \hbox{Sing} (\Psi)$, with $\hbox{Reg}(\Psi) \cap
\hbox{Sing}(\Psi) =0$, $\hbox{Reg} (\Psi)$ is a closed $d$-dimensional
submanifold of $M\backslash \hbox{Sing}(\Psi)$ and ${\cal
H}^d(\hbox{Sing}(\Psi))=0$, so one can define $[\Psi]=[\Psi\backslash
\hbox{Sing} (\Psi)]= [\hbox{Reg}(\Psi)]$. In the smooth category, it
is convenient to set up the following definition. Let $r\geq 1$ and
work in the ${\cal C}^r$ category, $r\geq 1$.

\medskip
\noindent
{\sc 1.1 Definition.} (\cite{HL}, \cite{SA}). A closed set $M$ in a
real manifold $X$ is called a ${\cal C}^r$-scarred manifold of
dimension $d$ if there exists a closed set $\sigma \subset M$ with
${\cal H}^d_{loc}(\sigma)=0$, such that $M\backslash \sigma$ is an
{\it oriented} ${\cal C}^r$-smooth $d$-dimensional submanifold of
$X\backslash \sigma$ with ${\cal H}_{loc}^d(M\backslash
\sigma)<\infty$.

\medskip
The smallest set $\sigma \subset M$ with the above properties is
called the {\it scarred set} of $M$. We adopt the notation 
$\sigma= \hbox{Sc}(M)$ and $\hbox{Reg}(M)=M\backslash \hbox{Sc}(M)$.
Nonetheless, if $M$ is ${\cal C}^r$-smooth, $d[M]=0$ of course does
not imply that $d[M\backslash \sigma]=0$ for a set $\sigma \subset M$
with ${\cal H}_{loc}^d(\sigma)=0$.

Let $M$ be a ${\cal C}^r$-scarred manifold of dimension $d$.  It
follows from Stokes's formula that, if ${\cal H}_{loc}^{d-1}(\hbox{Sc}(M))=0$,
then the current $[M]$ has no boundary, {\em i.e.} $d[M]=0$, in
particular if $M$ is a {\it complex} analytic set.  The current $[M]$
given by integration over $M\backslash \hbox{Sc}(M)$ is well defined,
but to retain the local behavior of a smooth current of integration,
one must add the condition that $d[M]=0$ locally, or globally, to have
a globally {\it closed} object, for example to solve a boundary
problem.

When $M$ is noncompact, the condition $d[M]=0$ shall mean the
following: $d([M]\cap U)=0$, for each oen set $U \subset \subset X$
with $\hbox{Int} \ \overline{U}=U$. One says that $d[M]=0$ locally.

\medskip
\noindent
{\sc 1.2. Definition.} $M$ is called a ${\cal C}^r$-scarred cycle if,
moreover, $d[M]=0$ locally.

\medskip
This condition is geometric in nature and is rather independent of the
measure-theoretic largeness that ${\cal
H}_{loc}^d(\hbox{Sc}(M))=0$. It corrects {\it globally} the
singularities (think of $\hbox{dim} \ M =1$).

\bigskip
\noindent
{\bf 2. Geometry of $M$ and CR currents.} Our purpose in this section
is to study the meaning of the notion of a CR meromorphic mapping
$(f,{\cal D}_f,[\Gamma_f])$ in the sense of Harvey and Lawson, in
particular the implications underlying that $[\Gamma_f]$ defines a
${\cal C}^r$-scarred manifold.  Following \cite{HL}, we begin by
establishing various useful equivalent formulations of the notion of
CR functions.  These definitions take place in the category of CR
objects and CR manifolds. Any locally embeddable CR manifold being
embeddable as a piece of a {\it generic} submanifold in $\C^{m+n}$,
{\em i.e.}  with $\hbox{CRdim} \ M =m$ and $\hbox{codim} \ M = n$, we
set up these concepts for $M$ {\it generic}.

Let $M$ be a ${\cal C}^r$-scarred CR manifold of type $(m,n)$ in
$\C^{m+n}$, {\em i.e.} of dimension $2m+n$, of CR dimension $m$ and of
codimension $n$. Denote by $t\in\C^{m+n}$ the coordinates on $\C^{m+n}$.
Near a point $p_0\in \hbox{Reg}(M)$, $M$ can be defined by cartesian
equations $\rho_j(t)=0$, $1\leq j\leq n$, where $\partial \rho_1
\wedge \cdots \wedge \partial \rho_n$ does not vanish on $M$. We then
have
$$T_p^cM=T_pM \cap JT_pM =\{X\in T_p M : \ \partial \rho_j(X)=0,
j=1,...,n\},$$ where $J$ denotes the usual complex structure on
$T\C^{m+n}$. Then $J$ can be extended to the complexification
$T_p^cM\otimes_{\R} \C$ with eigenvalues $\pm i$. Let $T_p^c M \otimes
\C=T^{1,0}M \oplus T^{0,1} M$ denote the decomposition into the
eigenspaces for $i$ and $-i$ respectively. Then there is a natural
$\C$-linear isomorphism from $T_p^cM$ to $T_p^{1,0}M$ given by the
correspondence $X \mapsto Z= \frac{1}{2} (X-iJX)$. Moreover, the
operation of complex conjugation is well-defined on $T_p^cM
\otimes_{\R} \C$ and we have $T_p^{1,0}M=\overline{T_p^{0,1}M}$.

Suppose now that $f: M\to \C$ is a function of class ${\cal C}^1$. $f$
is called a CR function if $\overline{L}f=0$, for every section
$\overline{L}$ of $T^{0,1}M$, {\em i.e.} $f$ is annihilated by the
antiholomorphic vectors tangent to $M$. Equivalently, the differential
$df$ is complex-linear at each point $p\in M$, $df(JX)=idf(X)$, for
all $X\in T_p^cM$.  The first definition continue to make sense for
the wider class of CR distributions on $M$.

To check a generalized definition in the distributional sense, let
$U\subset M$ be a small open set, let $l_1,...,l_m \in \Gamma(U,T^cM)$
and let $\lambda_1,...,\lambda_n \in \Gamma(U,TM)$ with
$l_1,Jl_1,...,l_m,Jl_m,\lambda_1,...,\lambda_n$ linearly independent.

These vector fields determine splittings $TU=T^cU \oplus \Lambda_U$
and $T^{*} U=T^cU^*\oplus \Lambda_U^*$ of the tangent bundle and the
cotangent bundle $T^*M$ restricted to $U$. The two spaces $T^cM$,
called the {\it complex tangent bundle}, and $H^0 M=(T^cM)^{\bot}$,
the annihilator of $T^cM$ in $T^*M$, called the {\it characteristic
bundle} of $M$, are canonical; the other two depend on the choice of a
splitting. Let
$l_1^*,Jl_1^*,...,l_m^*,Jl_m^*,\lambda_1^*,...,\lambda_n^*$ be the
dual covector fields. Naturally, if $f\in {\cal C}^1(U,\C)$ :
$$df=\sum_{j=1}^m(l_j(f)l_j^*+Jl_j(f)Jl_j^*)+ \sum_{k=1}^n
\lambda_k(f) \lambda_k^*.$$ Then one can define an induced
$\overline{\partial}$ operator on $M$ by
$$\overline{\partial}_M (f)=\sum_{j=1}^m \overline{L}_j(f)
\overline{L}_j^*,$$ where $\overline{L}_j=\frac{1}{2} (l_j +iJl_j)$
and $\overline{L}_j^*=(l_j^*-iJl_j^*)$, for $j=1,...,m$. Clearly, the
kernel of $\overline{\partial}_M$ is the ring of CR functions on $M$,
and the definition of $\overline{\partial}_M$ is independent of the
choice of local vector fields. However, the operator does depend on
the choice of the splitting of $TM$.

Note that if we extend the local vector fields used in the definition
of $\overline{\partial}_M$ above, to a neighborhood ${\cal U}$ of $U$
in $\C^{m+n}$, then we have $\overline{\partial}(f)= \sum_{j=1}^m
\overline{L}_j(f) \overline{L}_j^*+ \sum_{k=1}^n
\overline{\Lambda}_k(f)\overline{\Lambda}_k^*$, where
$\overline{\Lambda}_k=\frac{1}{2} (\lambda_k +iJ\lambda_k)$ and
$\overline{\Lambda}_k^*=(\lambda_k^*-iJ\lambda_k^*)$, for
$k=1,...,n$. If, furthermore $M=\{\rho_1=\cdots=\rho_n=0\}$ as above,
then along $M$ we can assume that
$\lambda_k^*=\partial\rho_k=(d\rho_k+id^c\rho_k)$.

\medskip
\noindent
{\sc 2.1. Proposition.} {\it Let $M$ be a piece of a ${\cal
C}^1$-smooth manifold with $\hbox{dim} \ M=2m+n$, 
closed in an open set ${\cal U}\subset
\C^{m+n}$. Then the following conditions are equivalent.

\hspace{1cm} {\rm (i)} $\hbox{dim}_{\C} T_pM\cap JT_pM=m$ for all
$p\in M$ {\rm ;}

\hspace{1cm} {\rm (ii)} $\int_M \alpha=0$ for all $(r,s)$-forms
$\alpha$ on ${\cal U}$ with $r+s=2m+n$ and $|r-s| > n$ {\rm ;}

\hspace{1cm} {\rm (iii)} $[M]=[M]_{m,m+n}+[M]_{m+1,m+n-1}+\cdots
+[M]_{m+n,m}$, where $[M]_{r,s}$ are the components of the current of
integration $[M]$ with respect to the Dolbeault decomposition{\rm ;}

\hspace{1cm} {\rm (iv)} $M$ is locally given by $n$ scalar equations
$x_j=h_j(y,w)$, $j=1,...,n$, in holomorphic coordinates $t=(w,z)$,
$w\in \C^m$, $z=x+iy\in \C^n$, with $h_j(0)=0$ and $dh_j(0)=0$.}

\medskip
The proof is omitted. When $M$ is ${\cal C}^1$-scarred, it is natural
to allow singularities also for maps defined over $M$.  The precise
formulation is due to Harvey and Lawson (\cite{HL}, II) and
favores the graph viewpoint. We transpose it in the CR category.

\medskip
\noindent
{\sc 2.2. Definition.} (\cite{HL}).  Let $M$ be a ${\cal C}^r$-scarred
submanifold of $\C^{m+n}$. Then a {\it ${\cal C}^r$-scarred} mapping
of $M$ into a complex manifold $Y$ is a ${\cal C}^r$-smooth map $f:
{\cal D}_f \to Y$ defined on an open dense subset ${\cal D}_f\subset
\hbox{Reg} (M)$ such that the closure $\Gamma_f$ of the graph
$\{(p,f(p))\in {\cal D}_f\times Y\}$ in $\C^{m+n} \times Y$ defines a
${\cal C}^r$-scarred {\it cycle} in $\C^{m+n} \times Y$, {\em i.e.}
$d[\Gamma_f]=0$.

\medskip
Scarred ${\cal C}^r$ mappings will constantly be denoted by $(f,{\cal
D}_f,[\Gamma_f])$, to remind precisely that they are not set-theoretic
maps.

\medskip
\noindent
{\sc 2.3. Definition.} (\cite{HL}, \cite{SA}).  $(f,{\cal
D}_f,[\Gamma_f])$ is called a ${\cal C}^r$-scarred {\it CR mapping}
if, moreover, $[\Gamma_f]$ is a ${\cal C}^r$-scarred CR cycle of
$\C^{m+n}\times Y$ of CR dimension $m$.

\medskip
One can go a step further in generalization.  Indeed, Harvey and
Lawson have introduced the notion of maximally complex
currents. Accordingly, {\it CR currents} arise as generalized currents
of integration on CR manifolds as follows.

\medskip
\noindent
{\sc 2.4. Definition.} Let ${\cal M}$ be a $d=(2m+n)$-dimensional
current with compact support on a $(m+n)$-dimensional complex manifold
$X$. ${\cal M}$ is called a generic current of type $(m,n)$ if the
Dolbeault components
$${\cal M}_{r,s}=0, \ \ \ \ \ r+s=2m+n, \ \ \hbox{for} \ \ |r-s|>n,$$
{\em i.e.} ${\cal M}(\alpha)=0$ for all $(r,s)$-forms $\alpha$ on $X$
where $|r-s| >n$.

\medskip
Let ${\cal M}$ be a closed generic current of type $(m,n)$ in an open
set ${\cal U} \subset \C^{m+n}$. Then ${\cal M}={\cal
M}^{0,n}+\cdots+{\cal M}^{n,0}$ and $d{\cal M}=0$ yield
$\overline{\partial} {\cal M}^{0,n}=0$, since $[d{\cal M}|^{0,n+1}=
\overline{\partial} {\cal M} ^{0,n}$, simply for reasons of
bidegree. Using this remark yields four equivalent definitions for a
${\cal C}^1$-smooth function to be CR. $[\Gamma_f]$ denotes the
current of integration over the closure of the graph of $f$. Since
$d[\Gamma_f]=0$, $\overline{\partial}[ \Gamma_f]^{0,n}=0$. The
variable $\zeta$ is used to denote a coordinate on $\C$, $f: M\to \C$.
Property (iv) below can be used as a new definition.
Let $\pi$ denote the projection $\C^{m+n}\times \C
\to \C^{m+n}$.

\medskip
\noindent
{\sc 2.5. Proposition.} (\cite{HL}). {\it Let $M$ be an oriented real
CR manifold of class ${\cal C}^1$ in an open set ${\cal U}\subset
\C^{m+n}$. Then, for any $f\in {\cal C}^1(M)$, the following
statements are equivalent:

\hspace{1cm} {\rm (i)} $f$ is a CR function on $M${\rm ;}

\hspace{1cm} {\rm (ii)} $\overline{\partial}_M(f)=0${\rm ;}

\hspace{1cm} {\rm (iii)} $\overline{\partial}(f[M]^{0,n})=0$, {\em
i.e.} $\int_{{\cal U}\cap M} f\overline{\partial}\varphi=0$, for all
$\varphi \in {\cal D}^{m+n,m-1}({\cal U})${\rm ;}

\hspace{1cm} {\rm (iv)}
$\overline{\partial}(\pi_*(\zeta[\Gamma_f]_{m+n,m}))= \pi_*(\zeta
\overline{\partial}[\Gamma_f]^{0,n})=0$.}

\medskip
{\it Proof.} Equivalence of (i) and (ii) is obvious. To prove that
(ii) implies (iii), it results from $\overline{\partial}_M(f)=0$ that
$f\overline{\partial}\varphi =\overline{\partial}(f\varphi)=
d(f\varphi)$, since $\partial(f\varphi)=0$ by bidegree considerations,
hence by Stokes' formula, $\int_Mf\overline{\partial}\varphi=\int_M
d(f\varphi)=0$. The converse is obtained by choosing adequate forms
$\varphi$.  To prove that (iii) is equivalent to (iv), notice that
$\overline{\partial}(f[M]^{0,n})=\overline{\partial}(\pi_*(\zeta
\overline{\partial}[\Gamma_f]^{0,n}))$, obviously.

The proof of Proposition 2.5 is complete.

\bigskip
\noindent
{\bf 3. CR meromorphic mappings.} The natural generalization of
meromorphy to CR category must include the appearance of indeterminacy
points, not only being smooth CR from $M$ generic to $\P^1(\C)$ or to
a projective algebraic manifold $Y$. The following definition was
devised by Harvey and Lawson and appears to be adequately large, but
sufficiently stringent to maintain the possibility of filling a
scarred maximally complex cycle with a holomorphic chain.

\medskip
\noindent
{\sc 3.1. Definition.} (\cite{HL}, \cite{SA}).  Let $M$ be a ${\cal
C}^r$-scarred generic submanifold of $\C^{m+n}$. Then a {\it CR
meromorphic mapping} is a ${\cal C}^r$-scarred CR mapping $(f,{\cal
D}_f,[\Gamma_f])$ with values in a projective manifold $Y$.

\medskip
{\it By definition}, a CR meromorphic mapping takes values
in a projective algebraic manifold.
In particular, if $Y=\P^1(\C)$, the closure $\Gamma_f$ of the graph of
$f$ over ${\cal D}_f$ defines a ${\cal C}^r$-scarred CR manifold of
type $(m,n)$ in $\C^{m+n} \times \P^1(\C)$ satisfying
$d[\Gamma_f]=0$. Since any projective $\P^k(\C)$ is birationaly
equivalent to a product of $k$ copies of $\P^1(\C)$ (\cite{HL}), we
can set $Y=\P^1(\C)$ without loss of generality.

\smallskip
{\it Remark.} We would like to mention that a map defined on a dense
open set ${\cal U} \subset \C^{m+n}$ with values in $\P^1(\C)$ is
meromorphic over ${\cal U}$ if and only if the closure $\Gamma_f$ of
its graph $\{(p,f(p)\in {\cal U} \times \P^1(\C)\}$ defines a ${\cal
C}^r$-scarred complex submanifold of ${\cal U}\times \P^1(\C)$.  This
justifies in a certain sense the above definition.

\smallskip 
Let $(t_1,...,t_{m+n},[\zeta_0:\zeta_1])=(t,\zeta)$ denote coordinates
on $\C^{m+n} \times \P^1(\C)$ and let $\pi: \C^{m+n} \times
\P^1(\C)\to \C^{m+n}$ denote the projection onto the first factor.

\medskip
\noindent
{\sc 3.2. Definition.} (\cite{SA}). A point $p\in M$ is called an {\it
indeterminacy point} if $\{p\}\times \P^1(\C) \subset \Gamma_f$. Denote by
$\Sigma_f=\{p\in M; \ \{p\}\times \P^1(\C) \subset \Gamma_f\}$ the
indeterminacy locus of $f$.

\medskip
The following two propositions are due to Sarkis (\cite{SA}).  The
first one is a clever remark about thinness of the indeterminacy set
$\Sigma_f$.  We expose his proof for completeness.

\medskip
\noindent
{\sc 3.3. Proposition.} ({\sc Sarkis}, \cite{SA}). {\it Let $M$ be a
${\cal C}^1$-scarred CR manifold of type $(m,n)$ in $\C^{m+n}$ and let
$(f,{\cal D}_f,[\Gamma_f])$ be a CR meromorphic mapping on $M$.
Then{\rm :}

\hspace{1cm} {\rm (i)} For almot all $a\in \P^1(\C)$, the level set
$\Lambda_a=\pi(\{\zeta=a\}\cap \Gamma_f)$ is a ${\cal C}^1$-scarred
$2$-codimensional submanifold of $M${\rm ;}

\hspace{1cm} {\rm (ii)} For every such $a$, the indeterminacy set
$\Sigma_f=\{p\in M; \ \{p\} \times \P^1(\C) \subset \Gamma_f\}$ of $f$
is a closed subset of $\Lambda_a$ with empty interior.}
\medskip

{\it Proof.} We begin by asserting that for almost all complex
$(m+n)$-dimensional linear subspaces $H$ of $\C^{m+n}\times \P^1(\C)$,
we have : 1. ${\cal H}^{d-2}(\Gamma_f \cap H) < \infty$ and 2.
$\Gamma_f^H:=\Gamma_f\cap H$ is a ${\cal C}^1$-scarred
$(2+n)$-codimensional real submanifold of $H$. This follows by known
facts from geometric measure theory, see \cite{HL}.  After a small
linear change of coordinates in $\C^{m+n}\times \P^1(\C)$, this holds
for almost every $a\in \P^1(\C)$ with $H=H_a=\{\zeta=a\}$. Write
$\Gamma_f^a=H_a\cap \Gamma_f$. Obviously, $\Gamma_f^a\subset M\times
\{a\}$ is a ${\cal C}^1$-scarred submanifold in $H_a$ if and only if
$\Gamma_f^a$ is a ${\cal C}^1$-scarred $2$-codimensional submanifold
of $M$.  This gives (i).

Assume by contradiction that $\Sigma_f$ contains a nonempty open set
${\cal L} \subset \hbox{Reg} (\Lambda_a)$, so ${\cal L} \times
\P^1(\C) \subset \Gamma_f$.  For reasons of dimension, ${\cal L}\times
\P^1(\C)\equiv \Gamma_f$ there. Indeed, $\hbox{dim}_{\R} ({\cal L}
\times \P^1(\C))=2+\hbox{dim}_{\R}{\cal L}=\hbox{dim}_{\R} \Gamma_f$.
Let $p_0\in {\cal L}$. That $\Gamma_f$ is vertical over ${\cal L}$
near $\hbox{Reg} (\Lambda_a)$ is impossible, since $\Gamma_f |_{{\cal
D}_f}$ is a ${\cal C}^1$-smooth graph over the dense open set ${\cal
D}_f \subset M$ whose closure contains $p_0$.

The proof of Proposition 3.3 is complete.

\smallskip
{\it Remark.} The small linear change of coordinates
above was necessary, since all the $H_a$ can be contained in the 
thin set of $H$ where 1 or 2 do not hold.

\smallskip
A classical observation is that to each pair consisting of volume form
$d\lambda_M$ on an oriented ${\cal C}^r$-scarred CR manifold $M$ and
an integrable function on $M$ is associated a distribution $T_f$ in a
natural way by $\left<T_f,\varphi\right>=\int_U f\varphi d\lambda_M$.
However, $T_f$ depends on $d\lambda_M$. There is associated the
transpose operator ${\:\!}^{\tau} \overline{L}$ of a CR vector field
$\overline{L}\in \Gamma(U,T^{0,1}M)$ with respect to $d\lambda_M$,
that is $\int_M \varphi \overline{L} (\psi) \ d\lambda_M=\int_M
{\:\!}^{\tau} \overline{L}(\varphi)\psi \ d\lambda_M$ for all
functions $\varphi, \psi$ with compact support.  Then $T_f$ is CR if
and only if $\left<T_f,{\:\!}^{\tau}\overline{L}(\varphi)\right>=0$ if
and only if $f$ is CR. A distribution $T$ on $M$ is called a CR
distribution if $\left<T,{\:\!}^{\tau}\overline{L}(\varphi)\right>=0$
for all $\varphi \in {\cal D}(M)$. Although
${\:\!}^{\tau}\overline{L}$ depends on the choice of $d\lambda_M$,
this annihilating condition is independent.  Indeed, given
$d\lambda_M^1$ and $d\lambda_M^2$, there always exists a function $a\in
{\cal C}^{\infty}(M, \C^*)$ with $d\lambda_M^2=ad\lambda_M^1$, so
${\:\!}^{\tau^2}\overline{L}(\varphi)= \frac{1}{a}
{\:\!}^{\tau^1}\overline{L} (a \varphi)$, whence the equivalence by
linearity of distributions.

The statement below and its proof are known if
$\hbox{Sc}(M)=\emptyset$, {\em i.e.} $f$ is ${\cal C}^1$; here, the
condition $d[\Gamma_f]$ helps in an essential way to keep it true in
the ${\cal C}^1$-scarred category.

\medskip
\noindent
{\sc Proposition 3.4.} ({\sc Sarkis}, \cite{SA}). {\it Let $M$ be a
${\cal C}^1$-scarred CR manifold of type $(m,n)$ in $\C^{m+n}$, let
$(f,{\cal D}_f,[\Gamma_f])$ be a CR meromorphic mapping on $M$ and let
$\Sigma_f=\{p\in M; \ \{p\}\times \P^1(\C) \subset \Gamma_f\}$. Then
there exists an order zero CR distribution $T_f$ on $M\backslash
\Sigma_f$ such that $T_f|_{{\cal D}_f} \equiv f$. In a chart $(U,\C)$
of $M\times \P^1(\C)$ with $(U\times \{\infty\})\cap
\Gamma_f|_U=\emptyset$, given a volume form $d\lambda_M$ on
$M\backslash \hbox{Sc}(M)$, $T_f$ is defined by
$$ 
[\Gamma_f](\zeta \pi_*(\varphi d\lambda_M))= \int_{\Gamma_f} \zeta
\pi^* (\varphi \ d\lambda_M),
$$
for all $\varphi\in {\cal C}_c^{\infty}(U)$.}

\medskip
{\it Proof.} As before, $\pi: U\times \C \to U$ denotes
$(z,\zeta)\mapsto z$. By assumption, $U\subset M\backslash
\Sigma_f$. Since $U\times \{\infty\} \cap \Gamma_f|_U = \emptyset$,
one has $\sup_{\zeta\in \Gamma_f|_U} |\zeta| < \infty$.  Let $V\subset
\subset U$ be open, $\overline{V}$ compact and let $\varphi\in {\cal
C}_c^{\infty}(V)$.  Then
$$|\left< T_f,\varphi\right>| = |\Gamma_f] (\zeta\pi^*(\varphi \
d\lambda_M))| \leq \sup_{\zeta\in \Gamma_f |_U} |\zeta| \ \ {\cal
H}^d(\Gamma_f \cap (V\times \P^1(\C)) ||\varphi||_{L^{\infty}(U)}$$
($d=\hbox{dim} M$), which proves that $T_f$ is a distribution of order
zero over $U$ ($T_f|_U\in L_{loc}^{\infty}(M)$.)

$T_f$ is clearly equal to the distribution associated with $f$ on the
open dense set ${\cal D}_f\subset M$ where $f$ is ${\cal
C}^1$. Indeed,
$$\forall \ \varphi\in {\cal C}_c^{\infty}(U), \ \ \ \ \ \left<T_f,
\varphi\right> = \int_{\Gamma_f \cap \pi^{-1}(U)} \zeta \pi^* (\varphi
d\lambda_M)= \int_U f\varphi d\lambda_M =\left<f,\varphi\right>.$$

Let now $\overline{L}\in \Gamma(U,T^{0,1}M)$ and complete the pair
$(L^*,\overline{L}^*)$ in a basis
$(L^*,\overline{L}^*,L_2^*,\overline{L}_2^*,...,$
$L_m^*,\overline{L}_m^*,\lambda_1^*,...,\lambda_n^*)$ of $T^*M$, so
that, furthermore,
$$d\lambda_M=\left(\frac{i}{2}\right)^m L^*\wedge \overline{L}^*
\wedge L_2^* \wedge \overline{L}_2^* \wedge \cdots \wedge L_m^* \wedge
\overline{L}_m^* \wedge \lambda_1^* \wedge \cdots \wedge
\lambda_n^*.$$ By Stokes' formula, $\int_U \varphi \overline{L} (\psi)
d\lambda_M= -\int_U \overline{L} (\varphi) \psi d\lambda_M$ for all
$\varphi, \psi \in {\cal C}^{\infty}_c(U)$, so that the transpose
${\:\!}^{\tau}\overline{L}$ equals $-\overline{L}$ in the above chosen
frame. One must prove that $\left< T_f, {\:\!}^{\tau} \overline{L}
(\varphi)\right>=0$ for all $\varphi \in {\cal C}^{\infty}_c(U)$. To
do so, notice that by introducing the $(2m+n-1)$-form $d\mu_M=
(\frac{i}{2})^m L^* \wedge L_2^* \wedge \overline{L}_2^* \wedge \cdots
\wedge \lambda_n^*$, one has
$$\overline{L} (\varphi) d\lambda_M =\overline{\partial}_M (\varphi
d\mu_M)= \overline{\partial}(\varphi d\mu_M)$$ on $M$, since
$\overline{\partial}_M=\sum_{j=1}^m \overline{L}_j(.)
\overline{L}^*_j$ and
$\overline{\partial}|_M=\overline{\partial}_M$. Therefore,
$$\left<T_f, \overline{L}(\varphi)\right>=[\Gamma_f] (\zeta \pi^*
(\overline{\partial}(\varphi d\mu_M)))= [\Gamma_f]^{0,n}
(\overline{\partial}(\zeta \pi^*(\varphi d\mu_M)))=0,$$ by the
above-noticed fact that $\overline{\partial}[\Gamma_f]^{0,n}=0$ and
since $\overline{\partial}(\zeta \pi^*(\varphi d\mu_M))$ is a
$(m+n,m)$-form on $U\times \C$.

The proof of Proposition 3.4 is complete.

\smallskip
{\it Remark.} In fact, $f$ induces an {\it intrinsic} CR current
$[C_f]$ on $M\backslash \Sigma_f$ by $[C_f](\alpha)= [\Gamma_f] (\zeta
\pi^* \alpha)$ in a chart as above.  CR distributions will be more
concrete for the properties of extendability.

\bigskip
\noindent
{\bf 4. Local extension of CR meromorphic mappings.} Let $\Omega$ be a
bounded domain with connected ${\cal C}^1$ boundary in $\C^n$, where
$n\geq 3$ and let $Y$ be a projective manifold. Harvey and Lawson
proved that any ${\cal C}^1$-scarred mapping $f: b\Omega \to Y$ which
satisfies the tangential Cauchy-Riemann equations at the regular
points of $f$ and such that $d[\Gamma_f]=0$ in the sense of currents
extends to a meromorphic map $F: \Omega \to Y$. By considering the
graph of $f$ over $b\Omega$, it is a corollary of the following
extension theorem.

\medskip
\noindent
{\sc Theorem.} ({\sc Harvey-Lawson}, \cite{HL}). {\it Let $(V,bV)$ be
a compact, complex, $p$-dimensional subvariety with boundary in
$\P^n(\C) \backslash \P^{n-q}(\C)$, where $bV$ is a scarred ${\cal
C}^1$-cycle whose regular points form a connected open set.  Then, if
$p>2q$, every scarred CR map of class ${\cal C}^1$ carrying $bV$ into
a projective manifold $Y$ extends to a meromorphic map $F: V \to Y$.}

\medskip
The case $\hbox{dim}_{\R} (b\Omega)=3$ and $b\Omega$ ${\cal C}^2$
follows from the work of Dolbeault-Henkin (\cite{DH}).

\medskip
\noindent
{\sc Theorem.} ({\sc Dolbeault-Henkin}, \cite{DH}).  {\it Let $\Omega$
be a bounded domain in $\C^2$, with $b\Omega$ of class ${\cal
C}^2$. Then every ${\cal C}^2$-smooth CR mapping $b\Omega \to
\P^1(\C)$ admits a meromorphic extension to $\Omega$.}

\medskip
In a forthcoming paper, Sarkis generalized the above result allowing
indeterminacies for $f$ CR meromorphic and a holomorphically convex
compact set $K$, in the spirit of Lupacciolu.

\medskip
\noindent
{\sc Theorem.} ({\sc Sarkis}, \cite{SA}). {\it Let $\Omega$ be a
relatively compact domain in a Stein manifold ${\cal M}$, $\hbox{dim}
\ {\cal M} \geq 2$, let $K= \widehat{K}_{{\cal H}({\cal M})}$ be a
holomorphically convex compact set and assume that $b\Omega \backslash
K$ is a connected ${\cal C}^1$-scarred hypersurface in ${\cal
M}\backslash K$. Then any CR meromorphic mapping $(f,{\cal D}_f,
[\Gamma_f])$ on $b\Omega \backslash K$ admits a unique meromorphic
extension to $\Omega \backslash K$.}

\medskip
We would like to mention that the above theorem is known for $f$ CR ${\cal
C}^1$ or CR meromorphic ${\cal C}^1$ without indeterminacies, by other
methods (\cite{PO2}), see Theorem 4.2 below.

\medskip
\noindent
{\it Global and local extension theorems.}  A general feature of
global extension of CR functions is that in many cases two independent
steps must be stated : I. Prove that $CR(M)$ extends holomorphically
(meromorphically) to a one-sided neighborhood ${\cal V}^b(M)$ (here,
$M$ is a ${\cal C}^1$-scarred hypersurface); II. Prove that the
envelope of holomorphy (meromorphy) of ${\cal V}^b(M)$ contains a
large open set, {\em e.g.} $\Omega$ if $M=b\Omega$. Step II is known
to be equivalent in both cases: the envelope of holomorphy and the
envelope of meromorphy of an open set coincide.

\medskip
\noindent
{\sc Theorem.} ({\sc Ivashkovitch}, \cite{IV}). {\it Let $Y$ be a
compact K\"ahler manifold and $f$ a meromorphic map from a domain
$\Omega$ in some Stein manifold into $Y$. Then $f$ extends to a
meromorphic map from the envelope of {\rm holomorphy}
$\widehat{\Omega}$ of $\Omega$ to $Y$.}

\medskip
Thus, every positive global extension theorem about CR functions
extends to be a result about meromorphic CR mappings, provided one can
prove by {\it local} techniques that they extend meromorphically to
open sets ${\cal V}^{b}(M)$ attached to real submanifolds $M\subset
\C^{m+1}$.  Indeed, the size of $\widehat{{\cal V}^b}(M)$ can be
studied by means of global techniques, {\em e.g.} integral formulas.
In most cases, including special results in partially convex-concave
manifolds, the disc envelope of such $M$ will contain some attached
open one-sided neighborhoods ${\cal V}^b(M)$ or wedges ${\cal W}$ with
edge $M$.

In this direction, a classical result is the Hartogs-Levi theorem:
{\it Let $\Omega \subset\subset \C^n$, $n\geq 2$, be a bounded domain
and let ${\cal V}(b\Omega)$ be an open neighborhood of its
boundary. Then holomorphic (meromorphic) functions on ${\cal
V}(b\Omega)$ extend holomorphically (meromorphically) to $\Omega$.}

Therefore, it is of great importance to answer the question of Henkin
and Sarkis (which was not raised by Harvey and Lawson in 1977): {\it Is there
a local version of the meromorphic extension phenomenon?} ({\em e.g.}
a Lewy extension phenomenon).  If the CR meromorphic mapping $(f,{\cal
D}_f,[\Gamma_f])$ does not possess indeterminacies, it is locally CR,
so the answer is positive.  We mention, however, that the most natural
notion of CR meromorphic maps it the one where indeterminacies really
occur, see Definitions 3.1 and 3.2.

Thus, a satisfactory understanding of CR meromorphicity involves the
local extension theory and various removable singularities theorems
(\cite{JO1}, \cite{JO3}, \cite{ME2}, \cite{PO1}, \cite{MP}).  This
paper is devoted to delineate some.

\medskip
\noindent
{\it CR meromorphy and removable singularities.} Let $M$ be a piece of
a generic submanifold of $\C^{m+n}$. The local holomorphic extension
phenomenon for CR(M) and ${\cal D}_{CR}'(M)$ as well arises at most
points of $M$, according to the theory of Tr\'epreau and Tumanov.

By a {\it wedge of edge} $M$ at $p_0\in M$, we mean an open set in
$\C^{m+n}$ of the form
$${\cal W}=\{z+\eta; \ z \in U, \eta \in C\},$$ for some open
neighborhood $U$ of $p_0$ in $M$ and some convex truncated open cone
$C$ in $T_{p_0}\C^{m+n}$, {\em i.e.} the intersection of a convex open
cone with a ball centered at $0$.

$M$ is called {\it minimal} at $p_0$ if the following property is
satisfied.

\medskip
\noindent
{\sc Theorem.} ({\sc Tr\'epreau}: $n=1$; {\sc Tumanov}: $n\geq 2$.)
{\it Assume $M\subset \C^{m+n}$ is generic, ${\cal C}^{2,\alpha}$
$(0<\alpha <1)$, $\hbox{CRdim} \ M =m \geq 1$, $\hbox{codim} \ M = n
\geq 1$ and let $p\in M$. Then there exists a wedge ${\cal W}_p$ of
edge $M$ at $p$ such that $CR(M)$, $L_{loc,CR}^1(M)$,
$L_{loc,CR}^{\infty}(M)$, ${\cal D}_{CR}'(M)$ extend holomorphically
to ${\cal W}_p$ if and only if there does not exist a CR manifold
$S\subset M$ with $S\ni p$ and $\hbox{CRdim} \ S= \hbox{CRdim} \ M$.}

\medskip
By Proposition 3.4, all components of CR meromorphic mappings
$(f,{\cal D}_f,[\Gamma_f])$ on $M$ behave locally like a CR
distribution outside the thin set $\Sigma_f$ of their indeterminacies,
therefore extendability properties hold everywhere outside $\Sigma_f$,
if $M$ is minimal at every point. Thus, to extend $f$ along steps I
and II, one is naturally led to the problem of propagating holomorphic
extension up to wedges over $\Sigma_f$. It appears that $\Sigma_f$ has
small enough size to be coverable by wedges. Namely, in the
hypersurface case, which has been intensively studied, all the
necessary results are already known : A wedge attached to $M\backslash
\Phi$, $\Phi\subset M$ closed, $\hbox{codim} \ M =1$, is simply an
open set ${\cal V}^b$ ($b=\pm$) containing at each point of
$M\backslash \Phi$ a one-sided neighborhood of $M$ such that
$\hbox{Int} \overline{{\cal V}^b}={\cal V}^b$. An open connected set
${\cal W}_0$ is called a {\it wedge attached to} $M\backslash \Phi$ if
there exists a continuous section $\eta: M \to T_M\C^{m+n}$ of the
normal bundle to $M$ and ${\cal W}_0$ contains a wedge ${\cal W}_p$ of
edge $M$ at $(p,\eta(p))$ for every $p\in M$.  A closed set
$\Phi\subset M$ is called {\it ${\cal W}$-removable} (${\cal
V}^b$-removable if $n=1$) if, given a wedge ${\cal W}_0$ attached to
$M\backslash \Phi$, there exists a wedge ${\cal W}$ attached to $M$
with holomorphic functions in ${\cal W}_0$ extending holomorphically
to ${\cal W}$.

J\"oricke in the ${\cal C}^2$-smooth case and then Chirka-Stout
weakening the smoothness assumption, using the profound solution by
Shcherbina of the three-dimensional Cauchy-Riemann Dirichlet problem
with continuous data, showed:

\medskip
\noindent
{\sc Theorem.} ({\sc J\"oricke}: ${\cal C}^2$ , \cite{JO1}; {\sc
Chirka-Stout}, \cite{CS}). {\it Let $M$ be a locally Lipschitz graphed
hypersurface in $\C^{m+1}$, let $\Sigma\subset M$ be a closed subset
with empty interior of a ${\cal C}^1$-scarred two-codimensional
submanifold $\Lambda \subset M$.  Then $\Sigma$ is ${\cal V}^b$
removable.}

\medskip
In the greater codimensional case also, to prove local extension of CR
meromorphic mappings one has in a natural way to prove ${\cal W}$
removability of $\Sigma_f$. Let us denote by ${\rm Sc}(\Sigma_f)$ the
scar set of a scarred manifold $\Lambda$ which contains 
the indeterminacy set by Section 3.

The main result of this paper is the following.

\medskip
\noindent
{\sc 4.1. Theorem.} {\it Let $M$ be a smooth
generic manifold in $\C^{m+n}$, $\hbox{CRdim} \ M =m
\geq 1$, $\hbox{codim} \ M =n \geq 1$ and assume that $M$ is minimal
at every point of $M$. Then there exists a wedge ${\cal W}_0$ attached
to $M$ such that all CR meromorphic mappings $(f,{\cal
D}_f,[\Gamma_f])$ extend meromorphically to ${\cal W}_0$ $(\Sigma_f$ is
${\cal W}$-removable$)$ under the following circumstances

{\rm (i)} $n=1$ $($hypersurface case$)$, $M$ is ${\cal C}^{2,\alpha}$ and
$($only$)$ globally minimal{\rm:}

{\rm (ii)} $M$ is ${\cal C}^{2,\alpha}$ and ${\cal H}^{d-3}({\rm Sc}(
\Sigma_f))=0${\rm ;}

{\rm (iii)} $M$ is ${\cal C}^{\omega}$ $($real analytic$)$.}

\medskip
{\it Remark.}  The wedge ${\cal W}_0$ is universal: it does not depend
on $(f,{\cal D}_f,[\Gamma_f])$.

\smallskip
{\it Remark.} The smoothness assumptions make Theorem 4.1 weaker in
the hypersurface case than the local meromorphic extension theorem
that follows from the theorem of J\"oricke-Chirka-Stout or than the
global theorem of Sarkis.  Nonetheless, $M$ need not be minimal at
every point, see Lemma 4.4 below.

\medskip
\noindent
{\it Applications: global meromorphic extension.} In the following
results, it is known that ${\cal V}^b(M)$ $M=b\Omega$, $M=b\Omega
\backslash \widehat{K}_{\overline{\Omega}}$, contain $\Omega$, $\Omega
\backslash \widehat{K}_{\overline{\Omega}}$ respectively (\cite{JO3},
\cite{PO1}). In the meromorphic case they were proved by
Sarkis (\cite{SA}, see also \cite{LU}, \cite{LT}, \cite{JO2},
\cite{PO1}).

\medskip
\noindent
{\sc 4.2. Theorem.} {\it Let $\Omega\subset \subset \C^{m+1}$ be a
${\cal C}^2$-bounded domain. Then any CR meromorphic mapping $(f,{\cal
D}_f,[\Gamma_f])$ on $b\Omega$ with values in $P^1(\C)$ extends
meromorphically to $\Omega$.}

\medskip
\noindent
{\sc 4.3. Theorem.} {\it Let $\Omega\subset \subset \C^2$ be a ${\cal
C}^2$-bounded domain and let $K \subset b\Omega$ be a compact
set. Then any CR meromorphic mapping on $b\Omega\backslash K$ with
values in $\P^1(\C)$ extends meromorphically to $\Omega \backslash
\widehat{K}_{\overline{\Omega}}$, where
$\widehat{K}_{\overline{\Omega}}=\{p\in \overline{\Omega}; \ |f(z)|
\leq \max_K |f|, \ \hbox{\rm for all functions } f\in {\cal H}({\cal
V}(\overline{\Omega}))\}$.}

\medskip
{\it Remark.} In the above two theorems, the hypersurface $M=b\Omega$
need not be everywhere minimal: CR(M) automatically extend
holomorphically to some ${\cal V}^b(M)$, since $M$ is known to be a
{\it single CR orbit} (\cite{JO2}).  To explain the phenomenon, we
need some definitions.

\smallskip

Let $M$ be a ${\cal C}^2$-smooth CR manifold.  The CR orbit of a point
$p\in M$ is the set of all endpoints of piecewise smooth integral
curves of $T^cM$ with origin $p$. CR orbits partition $M$. Sussmann
(see \cite{TR}, \cite{JO3}) showed that each CR orbit ${\cal O}_{CR}$
possesses a structure of a {\it smooth} ${\cal C}^1$ manifold making
the inclusion ${\cal O}_{CR} \to M$ an injective ${\cal C}^1$
immersion. By construction, each ${\cal O}_{CR}$ is a CR manifold with
$\hbox{CRdim} \ {\cal O}_{CR} = \hbox{CRdim} \ M$. Each CR manifold as
${\cal O}_{CR}$ is locally embeddable as a generic submanifold of some
$\C^N$, $N \leq m+n$. A CR manifold $M$ is called globally minimal if
$M$ consists of a single CR orbit.

The relevance of CR orbits to the extendability properties of CR
functions are owed to Tr\'epreau and yielded the following finest
possible {\it extension theorem}:

\medskip
\noindent
{\sc Theorem.} (\cite{TR},\cite{TU},\cite{JO2},\cite{ME1},\cite{PO1}).
{\it If $M$ is a globally minimal locally embeddable generic ${\cal
C}^{2,\alpha}$-smooth $(0<\alpha<1)$ manifold, there exists a wedge
${\cal W}_0$ attached to $M$ such that $CR(M)$, $L_{loc,CR}^1(M)$,
$L_{loc,CR}^{\infty}(M)$, ${\cal D}_{CR}'(M)$ extend holomorphically
to ${\cal W}_0$.}

\medskip
{\it Proof of Theorems 4.2 and 4.3.} It is a fact that in the
hypersurface case, thin sets as $\Sigma_f$ do not perturb CR orbits:

\medskip
\noindent
{\sc 4.4. Lemma.} {\it Let $M$ be a ${\cal C}^2$ hypersurface in
$\C^{m+1}$ and let $\Sigma$ be a closed subset with nonempty interior
of some ${\cal C}^1$-scarred two-codimensional submanifold $\Lambda
\subset M$.  Then, for all CR orbits ${\cal O}_{CR}\subset M$, ${\cal
O}_{CR} \backslash ({\cal O}_{CR} \cap \Sigma)$ is a single CR orbit
of $M\backslash \Sigma$.}

\medskip
{\it Proof.} The real dimension of a ${\cal O}_{CR}$ is $\geq 2m$ and
$\leq 2m+n=2m+1$ if $n=1$. So $\Sigma$ is too small to make
obstruction to an orbit.  However, the lemma can fail in codimension
$\geq 2$.

\medskip
Thus, Theorems 4.2 and 4.3 rely on the following properties.

\medskip
\noindent
{\sc Proposition.} {\it Let $M=b\Omega$ or $b\Omega \backslash
\widehat{K}_{\overline{\Omega}}$, $\hbox{codim} \ M =1$, $M$ ${\cal
C}^2$ and let $(f,{\cal D}_f,[\Gamma_f])$ be a CR meromorphic mapping
on $M$. Then

\hspace{1cm} {\rm (i)} $M$ is a single CR orbit $($\cite{JOMZ}$)${\rm
;}

\hspace{1cm} {\rm (ii)} $M\backslash \Sigma_f$ is a single CR orbit,
hence $f$ extends meromorphically to ${\cal V}^b(M\backslash
\Sigma_f)${\rm ;}

\hspace{1cm} {\rm (iii)} $\Sigma_f$ is ${\cal V}^b$ removable, hence
$f$ extends meromorphically to ${\cal V}^b(M)$.}

\medskip
One then concludes 4.2 and 4.3 with the Ivashkovitch theorem.

\medskip
\noindent
{\it ${\cal W}$-removability.} Theorem 4.1 is reduced to the ${\cal
W}$-removability of $\Sigma_f$. By Proposition 3.3, $\Sigma_f$ is a
closed subset with empty interior of some ${\cal C}^1$-scarred
two-codimensional submanifold $\Lambda \subset M$, $\Sigma_f \subset
\hbox{Sc}(\Lambda) \cup \hbox{Reg}(\Lambda)$. Write $\Sigma_f=E\cup
\Phi$, $\Phi= \hbox{Reg} (\Lambda) \cap \Sigma_f$, $E=
\hbox{Sc}(\Lambda) \cap \Sigma_f$, ${\cal H}^{d-2}(E)=0$. $\Phi$ is
already known to be removable.

\medskip
\noindent
{\sc Theorem.} (\cite{ME2}, \cite{MP}).  {\it Let $M$ be a ${\cal
C}^{2,\alpha}$-smooth $(0<\alpha<1)$ generic manifold in $\C^{m+n}$,
minimal at every point, $\hbox{CRdim} \ M =m \geq 1$, and let
$N\subset M$ be a connected ${\cal C}^1$-smooth submanifold with
$\hbox{codim}_M N= 2$. Then every proper closed subset $\Phi \subset
N$ is ${\cal W}$-removable.}

\medskip
The purpose of Section 5 is to establish:

\medskip
\noindent
{\sc 4.5. Theorem.} {\it Let $M$ be a ${\cal C}^{2,\alpha}$-smooth
generic manifold in $\C^{m+n}$, globally minimal
with $\hbox{CRdim} \ M =m \geq 1$. Then

{\rm (i)} If $n=1$ or $M$ is ${\cal C}^{\omega}$, then any closed 
$E\subset M$ with 
${\cal H}_{loc}^{2m+n-2}(E)=0$ is ${\cal W}$ removable{\rm ;}

{\rm (ii)} If $M$ is minimal at every point, then any
closed subset $E\subset M$ with 
${\cal H}_{loc}^{2m+n-3}(E)=0$ is ${\cal W}$ removable.}

\medskip
{\it Remark.} Dinh and Sarkis obtained Theorem 4.5 assuming that $M$
is of type one in the sense of Bloom-Graham, {\em i.e.}  the first
order Lie brackets of vector fields in $T^cM$ generate $TM$,
$[T^cM,T^cM]=TM$, for $M$ ${\cal C}^4$-smooth (\cite{DS}).

\medskip
\noindent
{\it $L^{\rm p}$ removability.} Let $M$ be a locally embeddable CR
manifold of class ${\cal C}^2$.  A closed subset $\Phi$ of $M$ is
called {\it $L^{\rm p}$ removable}, $p\geq 1$, if each function $f\in
L^{\rm p}_{loc}(M)$ which satisfies the Cauchy-Riemann equations
$Lf=0$ (in the distribution sense) on $M\backslash \Phi$ satisfies the
equation $Lf=0$ on the whole of $M$, or, for short, if
$$L_{loc,CR}^{\rm p}(M\backslash \Phi) \cap L_{loc}^{\rm p} (M)=
L_{loc, CR}^{\rm p}(M).$$ The authors have proved in \cite{MP} that
$L^{\rm p}$ removability holds if ${\cal W}$ removability holds, for
closed subsets $\Phi \subset M$ with ${\cal H}_{loc}^{d-2}(\Phi) <
\infty$. Therefore:

\medskip
\noindent
{\sc 4.6. Theorem.} {\it Theorems 4.5 and 5.1 are true for $L^{\rm p}$
removability, $1\leq {\rm p} \leq \infty$.}

\bigskip
\noindent
{\bf 5. Removable singularities.}  As pointed out in Section 4, the
relationships between extendability properties of CR functions and the
geometry of a CR manifold are adequately reflected by its CR
orbits. It is thus natural to state a removable singularities theorem
in the most general context.

\medskip
\noindent
{\sc 5.1. Theorem.} {\it Let $M$ be a smooth
generic globally minimal manifold in $\C^{m+n}$ with
$\hbox{CRdim} \ M = m \geq 1$, $\hbox{codim} \ M = n \geq 1$ and
$\hbox{dim} \ M=d=2m+n$. 
Then every closed subset $E$ of $M$ such that
$M\backslash E$ is globally minimal is
${\cal W}$ removable under each of the following conditions

{\rm (i)} $n=1$, $M$ is ${\cal C}^{2,\alpha}$ and ${\cal H}^{d-2}(E)=0${\rm ;}

{\rm (ii)} $M$ is ${\cal C}^{2,\alpha}$ and ${\cal H}^{d-3}(E)=0${\rm ;}

{\rm (iii)} $M$ is ${\cal C}^{\omega}$.}

\medskip
{\it Proof.} Following the scheme of proof devised in \cite{ME2} and
\cite{MP}, we present the developement of the proof of (i), (ii) and (iii)
in five essential steps. Let $E\subset M$ be closed with 
${\cal H}^{d-2}(E)=0$.

\smallskip
\noindent
{\it Step one: Reduction to the removal of a point.} By assumption,
$M\backslash E$ is globally minimal. Then, according to the extension
theorem, CR function are wedge extendable at every point of
$M\backslash E$. Nonetheless, the direction of the above mentioned
wedges can suffer of discontinuities. Fortunately, {\it the edge of
the wedge theorem} enables one to fill in larger wedges by means of
attached analytic discs at points of discontinuity. Therefore, there
exists a wedge ${\cal W}_0$ attached to $M\backslash E$ to which
$CR(M)$ holomorphically extends.

Using a ${\cal C}^{2,\alpha}$-smooth partition of unity on
$M\backslash E$, we can deform $M$ inside ${\cal W}_0$ over
$M\backslash E$ in a ${\cal C}^{2,\alpha}$-smooth manifold
$M^d$. Then, instead of a function $f\in CR(M\backslash E)$, we get a
function $f$, holomorphic in a neighborhood $\omega$ ($\equiv {\cal
W}_0$) of $M^d\backslash E$ in $\C^{m+n}$.  The aim will be
subsequently to prove that such holomorphic functions extend into a
wedge ${\cal W}_1^d$ attached to $M^d$. The construction will depend
smoothly on $d$, so that letting $d$ tend to zero, one obtains a wedge
${\cal W}_1$ attached to $M$ (for details, see Section 5 in
\cite{ME2}).

The first key point is that the {\it continuity principle} along
analytic discs with boundaries in $\omega$ can now be exploited to
show that the envelope of holomorphy of $\omega$ contains a wedge
${\cal W}_1^d$ attached to $M^d$.

Let $\Delta$ denote the unit disc in $\C$ and $b\Delta$ its boundary,
the unit circle.  An embedded analytic disc $A$ attached to $M$ is
said to be {\it analytically isotopic} to a point in $M$ if there
exists a $C^1$-smooth mapping $(s,\zeta)\mapsto A_s(\zeta), 0\leq s
\leq 1, \zeta \in \overline{\Delta}$, such that $A_0=A$, each $A_s$ is
an embedded analytic disc attached to $M$ for $0\leq s <1$ and $A_1$
is a constant mapping $\overline{\Delta} \to \{pt\} \in M$.  Using
Cauchy estimates and controlling connectedness, it is possible to
prove (the embedding condition yields monodromy, \cite{ME2},
Proposition 3.2):

\medskip
\noindent
{\sc 5.2. Proposition.}  {\it Let $M$ be generic, $C^{2,\alpha}$, let
$\Phi$ be a proper closed subset of $M$ and let $\omega$ be a
neighborhood of $M\backslash \Phi$ in $\C^{m+n}$. If an embedded disc
$A$ attached to $M\backslash \Phi$ is analytically isotopic to a point
in $M\backslash \Phi$, then, there exists a neighborhood ${\cal
V}(A(\overline{\Delta}))$ in $\C^{m+n}$ such that, for each function
$f\in {\cal H}(\omega)$, there exists a function $F\in {\cal H}({\cal
V}(A(\overline{\Delta})))$ such that $F=f$ in a neighborhood of
$A(b\Delta)$.}

\medskip
Call a point $p\in E$ ${\cal W}$ removable if there exists a wedge
${\cal W}_p$ of edge $M^d$ at $p$ with ${\cal H}(\omega)$ extending
holomorphically to ${\cal W}_p$.

Define
$${\cal A}= \{\Psi \subset E \ \hbox{closed} ; M\backslash \Psi \
\hbox{is globally minimal and} \ M\backslash \Psi \ \hbox{is} \ {\cal
W} \ \hbox{removable in} \ M\backslash \Psi \ \}$$ and define
$E_{nr}=\cap_{\Psi\in {\cal A}} \Psi$, the non removable part of $E$.
Then $M\backslash E_{nr}$ is globally minimal too.  By deforming $M^d$
in a manifold $(M^d)^{d_1}$ over $E\backslash E_{nr}$, we can assume
that we must remove $E_{nr}$ for ${\cal H}({\cal
V}((M^d)^{d_1}\backslash E_{nr}))$ instead of $E$. But $E_{nr}$ is the
smallest non-removable subset of $E$ keeping $(M^d)^{d_1}\backslash
E_{nr}$ globally minimal.  Assuming that $E_{nr}\neq \emptyset$, we
shall now reach a contradiction by showing that a point $p_1 \in
E_{nr}$ is ${\cal W}$ removable.  We take now the notations $E$ and
$M$ instead of $(M^d)^{d_1}$. Thus, to prove Theorem 5.1, it is
sufficient to prove that the new $E$ is removable near one of its
points.

According to Lemma 2.3 in \cite{MP}, the fact that $M\backslash E$ is
globally minimal and the existence of chains of infinitesimally small
analytic discs approximating integral curves of $T^cM$ insure the
existence of a generic manifold $M_1$ of codimension one in $M$
through a point $p_1 \in E$ such that $T_{p_1} M_1\not\supset
T_{p_1}^cM$ and $E \subset M_1^-$, the closed negative side of $M_1$
in $M$, near $p_1$. Let us quote this (elementary) differential
geometric statement as: {\it Let $M$ be a ${\cal C}^2$ manifold, let
$K\subset TM$ be a ${\cal C}^1$ subbundle, let $E\subset M$ be a
closed nonempty set and assume that $M$ and $M\backslash E$ are both
single $K$-orbits. Then there exists a point $p_1\in E$ and
$M_1\subset M$ a ${\cal C}^1$ hypersurface with $p_1\in M_1$,
$T_{p_1}M_1 \not\supset K(p_1)$ and $E\subset M_1^-$ near $p_1$.}

Finally, by the definition of ${\cal A}$, and by disposition of $M_1$,
$E\subset M_1^-$, it suffices to show that $p_1$ is ${\cal W}$
removable.  Indeed, for a small neighborhood ${\cal V}(p_1)$ of $p_1$
in $M$, $(M\backslash E) \cup {\cal V}(p_1)$ is globally minimal, as
$T_{p_1} M_1 \not\supset T_{p_1}^cM$.  Thus, to prove Theorem 5.1, it
is sufficient to prove that a neighborhood of $p_1\in E$ is ${\cal
W}$-removable.

\medskip
\noindent
{\it Step two: Existence of a disc.} Let $p_1\in E$ as above and
choose holomorphic coordinates $(w,z)=(w_1,...,w_m,z_1,...,z_n)$,
$z=x+iy$ on $\C^{m+n}$ such that $p_1=0$, $T_0M=\{x=0\}$,
$T_0^cM=\{z=0\}$, $M$ is given by $n$ scalar equations $x=h(y,w)$, in
vectorial notation, $h(0)=0$, $dh(0)=0$ and $M_1$ is given in $M$ by
the supplementary equation $u_1=k(v_1,w_2,...,w_m,y)$, for a ${\cal
C}^{2}$-smooth $k$ with $k(0)=dk(0)=0$.  We denote $M_1^-=\{u_1 \leq
k(v_1,w_2,...,w_m,y)\}$.

Our first construction of analytic discs attached to $M$ proceeds as
follows.

\medskip
\noindent
{\sc 5.3. Lemma.} (\cite{MP}, Lemma 2.4).  {\it There exists an
embedded analytic disc $A\in C^{2,\beta}(\overline{\Delta})$ with
$A(1)=p_1$, $A(b\Delta)\backslash \{1\}\subset M\backslash M_1^{-}$
and $\frac{d}{d\theta}|_{\theta=0} A(e^{i\theta})= v_0\in
T_{p_1}M_1$.}

\medskip
It suffices to take, for small $\rho_1 >0$, the disc $A(\zeta)$ with
$W_{\rho_1}(\zeta)= (\rho_1(1-\zeta),0,...,0)$ and with $Y$ component
satisfying Bishop's equation $Y_{\rho_1}=T_1h(Y_{\rho_1},W_{\rho_1})$
on $b\Delta$ (here, $T_1$ denotes the Hilbert transform $L^2(b\Delta)
\to L^2(b\Delta)$ vanishing at 1, $(T_1u)(1)=0$).

\smallskip

Therefore, removability of $p_1$ wille be a consequence of the following
proposition 5.4, proved below, and of  Theorem 5.10 below. 
This is the main technical part
of the article. This proposition provides extension outside a
thin set ${\cal E}_{\Phi_E}$ which is studied below and which
lives in an {\it open} (wedge) set ${\cal W}$ of $\C^{m+n}$.

\medskip
\noindent
{\sc 5.4. Proposition.} {\it Let $M$ be generic,
$C^{2,\alpha}$-smooth, let $p_1\in M$, let $E\subset M$ be a closed
subset with ${\cal H}_{loc}^{d-2}(E)=0$, let $p_1\in E$ and assume that
there exists a one codimensional generic ${\cal C}^2$-smooth manifold
$M_1 \subset M$ such that $E \subset M_1^-$ and let $\omega$ be a
neighborhood of $M \backslash E$ in $\C^{m+n}$.  Let $A\in
C^{2,\beta}(\overline{\Delta})$ be a sufficiently small embedded
analytic disc attached to $M$, $A(1)=p_1$,
$\frac{d}{d\theta}|_{\theta=0} A(e^{i\theta})= v_0\in T_{p_1}M_1$,
with $A(b \Delta \backslash \{1\}) \subset M \backslash M_1^{-}$.
Then for each $\varepsilon >0$, there exist $v_{00}\in T_{p_1}M_1$
with $|v_{00}-v_0|<\varepsilon$, $v_{00}\not\in T_{p_1}^cM$, a
wedge ${\cal W}$ of edge $M$ at $(p_1, Jv_{00})$ 
and a closed set ${\cal E}_{\Phi_E}$ which is ${\cal C}^{2,\alpha}$
foliated by complex curves with ${\cal H}^{2m+2n-1}({\cal E}_{\Phi_E})=0$
such that for every
holomorphic function $f \in {\cal H}(\omega)$ there exists a function
$F \in {\cal H}(\omega\cup ({\cal W}\backslash {\cal E}_{\Phi_E}))$ 
with $F=f$ in the intersection of ${\cal
W}\backslash {\cal E}_{\Phi_E}$ 
with a neighborhood of $M\backslash E$ in $\C^{m+n}$.}

\medskip
{\it Remark.} For any $e\geq 2$, we obtain statement 5.4 above
for ${\cal H}^{d-e}(E)=0$
with ${\cal H}^{2m+2n-e+1}({\cal E}_{\Phi_E})=0$ 
(${\cal H}^{d-2}(E)=0$ is crucial for isotopies, see 5.8 below).

\medskip
\noindent
{\it Step three: Maximal families of analytic discs.}  This step
consists in including the above analytic disc in a very large
parameterized family of analytic discs obtained by varying the $W$
component, and its approximate radius, the base point $A(1)=p$ in a
small neighborhood of 0, and the point $A(-1)$ in $\omega$
(\cite{ME2}, \cite{MP}).

Let $\mu=\mu(y,w)$ be a $C^{\infty },$ $\R$-valued function with
support near the point $(y(-1), w(-1))$ that equals 1 there and let
$\kappa: \R^{n} \to \R^{n}$ be a $C^{\infty }$ function with
$\kappa(0)=0$ and $\kappa'(0)=Id$.  We can assume that the supports of
$\mu$ and $\kappa$ are sufficiently concentrated in order that every
manifold $M_t$ with equation \begin{equation}
x=H(y,w,t)=h(y,w)+\kappa(t)\mu(y,w) \end{equation} is contained in
$\omega$ and the deformation is localized in a small neighborhood of
$A(-1)$ in $\C^{m+n}$.  Let $\chi= \chi(\zeta)$ be a smooth function
on the unit circle supported in a small neighborhood of $\zeta=-1$.

We consider the disc with $W$ component
$$W_{\tau,a,\rho,p}(\zeta)=(e^{i\tau}(\rho_1-\rho\zeta)+w^0_1,
(1-\zeta)a_1\rho/\rho_1+w_2^0,...,
(1-\zeta)a_{m-1}\rho/\rho_1+w_m^0),$$ where $a\in \C^{m-1}$ runs in a
small neighborhood ${\cal A}$ of 0, $0\leq \rho \leq \rho_1$, $w^0\in
\C^m$, a point $p\in M$ runs in a neighborhood of $0$ is represented
by its coordinates $(w^0,y^0)$ and with $Y$ component
$Y_{t,\tau,a,\rho,p}$ which is the solution of Bishop's equation with
parameters
$$Y_{t,\tau,a,\rho,p}=T_1H(Y_{t,\tau,a,\rho,p},
e^{i\tau}(\rho_1-\rho\zeta)+w^0_1,(1-\zeta)a_1\rho/\rho_1
+w_2^0,..., (1-\zeta)a_{m-1}\rho/\rho_1+w_m^0,t\chi)+y^0,$$ which
exists and depends in a ${\cal C}^{2,\beta}$-smooth fashion on
$(t,\tau,a,\rho,p,\zeta)$, for all $0<\beta<\alpha$.  Then
$A_{t,\tau,a,\rho,p}(1)=p$. When $\tau=0$, $a=0$, $\rho=\rho_1$ and
$p=0$, simply denote $A_{t,0,0,\rho_1,0}$ by $A_t$.

Let us recall that the normal deformations of $A$ near $A(-1)$ in
$\omega$ can be chosen in order that the inner tangential direction
$-\partial A_t /\partial \zeta (1)$ will describe a {\it whole} open
cone in the normal bundle to $M$ at $A(1)=0$.

Let $\Pi$ denote the canonical bundle epimorphism $\Pi: T\C^{m+n} |_M
\to T\C^{m+n} |_M / TM$ and consider the $C^{1,\beta}$ mapping
\begin{equation} D: \ \ \ \ \ \R^{n} \ni t \ \longmapsto \ \Pi \left(
- \frac{\partial A_t}{\partial \zeta} (1) \right) \in T_0 \C^{m+n} /
T_0 M \simeq \R^{n}.  \end{equation} We refer to \cite{MP} for a proof
of the following.

\medskip
\noindent
{\sc 5.5. Lemma.} ({\sc Tumanov} \cite{TU}). {\it $\chi$ can be chosen
such that $ \ rk \ D'(0)= n$.}

\medskip
This statement is more or less equivalent to the fact that the union
of the discs describes a wedge of edge $M$ at $0$. We also have:

\medskip
\noindent
{\sc 5.6. Lemma.} (\cite{MP}). {\it $\chi$ can be chosen such that the
following holds: there exist $\tau_0 >0$, ${\cal T}$ a neighborhood of
$0$ in $\R^n$ and ${\cal A}$ a neighborhood of $0$ in $\C^{m-1}$ such
that the set \begin{equation} \Gamma_0= \{s
\frac{dA_{t,\tau,\rho_1,a,0}}{d\theta} (1); \ s>0, t \in {\cal T},
\tau\in I_{\tau_0}, a\in {\cal A}\} \end{equation} is a
$(2m+n)$-dimensional open connected cone with vertex $0$ in $T_0M$.}

\medskip

For convenience, we shall allow us to shrink any open neighborhoods
arising in the next constructions without explicit mention. By reasons
of rank, the geometric meaning will be clear for sufficiently small
parameters.

\medskip
\noindent
{\it Step four: Isotopies.} The main hypothesis so far is
${\cal H}^{d-2}(E)=0$. Boundaries of analytic discs
$A_{t,\tau,a,\rho,p}(b\Delta)$ are embedded ${\cal
C}^{2,\beta}$-smooth copies of $S^1$ in $M$, so one expects naturally
that $A_{t,\tau,a,\rho,p}(b\Delta)\cap E = \emptyset$
generically. Furthermore, an isotopy property is required as stated in
Proposition 5.2, in order to have that ${\cal H}(\omega)$ extends
holomorphically to ${\cal V}(A_{t,\tau,a,\rho,p}(\overline{\Delta}))$.

To prove an isotopy lemma, we recall briefly some facts concerning
Hausdorff measures, taken from the very clear exposition of Chirka
\cite{CH}.

Let $E$ be an arbitrary subset of a metric topological space. For
$\delta>0$, define
$${\cal H}_{\delta}^s(E)=\inf \{\sum_{j=1}^{\infty} r_j^s; \ E \
\hbox{is covered by} \ \bigcup_{j=1}^{\infty} B_j, B_j= \ \hbox{balls
of radius} \ r_j \leq \delta\}.$$ Clearly, ${\cal H}_{\delta}^s (E)
\leq {\cal H}_{\delta'}^{s}(E)$, for $\delta ' \leq \delta$, so the
limit ${\cal H}^s(E)=\lim_{\delta\to 0} {\cal H}_{\delta}^s(E)$ exists
in $[0,\infty]$ and is called the $s$-dimensional Hausdorff measure of
$E$. The important property is that there exists a critical $\gamma
\geq 0$, the so-called Hausdorff dimension of $E$, such that ${\cal
H}^s(E)=\infty$ for all $s < \gamma$ and ${\cal H}^{s}(E)=0$ for all
$s> \gamma$, and the value of ${\cal H}^{\gamma}(E)$, if in
$(0,\infty)$, is not important.

This notion of dimension especially applies in the category of ${\cal
C}^1$ manifolds. Let $M$ and $N$ be connected real riemannian
manifolds, of class ${\cal C}^1$, $\hbox{dim} \ M =d \geq 1$ and let
$E \subset M$ be closed.

\medskip
\noindent
{\sc 5.7. Proposition.} (\cite{CH}, 346-352).  {\it {\rm (i)} ${\cal
H}^0(E)= \hbox{Card}(E)${\rm ;}

\hspace{1cm} {\rm (ii)} ${\cal H}^d$ coincides with the outer Lebesgue
measure on $M${\rm ;}

\hspace{1cm} {\rm (iii)} If ${\cal H}^{d-1}(E)=0$, then $M\backslash
E$ is locally connected{\rm ;}

\hspace{1cm} {\rm (iv)} Let $\pi: M \to N$ be a ${\cal C}^1$-smooth
map and let $E \subset M$ such that ${\cal H}^s(E)=0$, for a $s \geq
e=\hbox{dim}\ N$. Then ${\cal H}^{s-e}(E\cap \pi^{-1}(y))=0$ for
$d\lambda_N$-almost all $y\in N$.}

\medskip
Properties (i), (iii) and (iv) are naturally involved in the proof of
the following.

\medskip
\noindent
{\sc 5.8. Lemma.} {\it Let $E \subset M_1^-$ be a closed set with
${\cal H}_{loc}^{d-2} (E)=0$. Then for all small $(t,\tau,a,\rho,p)$,
each disc with $A_{t,\tau,a,\rho,p}(b\Delta) \cap E = \emptyset$ is
analytically isotopic to a point in $M\backslash E$.}

\medskip
{\it Proof.} During this proof $t$, $\tau$, $a$ and $w_1^0$ are fixed.
Then there exist $0< \rho_1$, $I_{\rho_1}=(0,\rho_1)$, a neighborhood
${\cal V}^*$ of $0$ in $\C_{w^*}^{m-1}$, $w^*=(w_2,...,w_m)$ and a
neighborhood $\Upsilon$ of $0$ in $\R^n$ such that the mapping (note
that $p$ is parameterized by $(w_1^0,w^{0*},y^0))$:
$${\cal S} \times b\Delta = I_{\rho_1} \times {\cal V}^* \times
b\Delta \ni (\rho, w^{0*}, y^0, \zeta) \mapsto
A_{t,\tau,a,\rho,p}(\zeta)\in M$$ is an embedding. Indeed, this
follows by differentiating Bishop's equation, noting first that $\partial
Y_{0,0,0,0,0}/\partial y^0=\hbox{Id}$, $\partial W_{0,0,0,0} /\partial
w^{0*} =\hbox{Id}$, that $I_{\rho_1} \times b\Delta \ni (\rho,\zeta)
\mapsto (\rho_1-\rho \zeta)\in \C$ is an embedding and recognizing
that $A_{t,\tau,a,\rho,p}(\zeta)$ is ${\cal C}^{2,\beta}$ with respect
to all variables. This exhibits a foliation of an open set in $M$ by
${\cal C}^{2,\beta}$ {\it real} discs ${\cal D}_{t,\tau,a,p}={\cal
D}_{t,\tau,a,w_1^0,s}$, $s=(w^{0*}, y^0)$, where
$${\cal D}_{t,\tau,a,w_1^0,s}=\{A_{t,\tau,a,\rho,p}(\zeta)\in M; \
0\leq \rho < \rho_2, \zeta \in b\Delta\}.$$ Now, since ${\cal
H}^{d-2}(E)=0$, the set
$${\cal S}\backslash {\cal S}_{E,t,\tau,a,w_1^0}=\{s\in {\cal S}; \
{\cal H}^0 ({\cal D}_{t,\tau,a,w_1^0,s}\cap E)=0\}$$ is a full measure
$(d-2)$-dimensional subset of ${\cal S}=I_{\rho_1} \times {\cal V}^*
\times {\cal Y} \simeq \R^{d-2}$, by Proposition 5.7 (iv).  This shows
that ${\cal D}_{t,\tau,a,w_1^0,s}\cap E=\emptyset$ for $d\lambda_{\cal
S}$-almost $s=(w^{0*}, y^0)\in {\cal S}$, where $t, \tau, a$ and
$w_1^0$ are fixed. Clearly, the mapping
$$I_{\rho_2} \cup \{0\} \times \overline{\Delta} \ni (\rho, \zeta)
\mapsto A_{t,\tau,a,\rho,p}(\zeta)\in \C^{m+n}$$ yields an analytic
isotopy of the analytic discs $A_{t,\tau,a,\rho,p}(\zeta)$, for all
$0<\rho \leq \rho_1$, provided $s\in {\cal S} \backslash {\cal
S}_{E,t,\tau,a,w_1^0}$.  It remains to show that discs such that
$A_{t,\tau,a,\rho,p}(b\Delta) \cap E =\emptyset$ but
$${\cal D}_{t,\tau,a,w_1^0,s} \cap E \neq \emptyset$$ are also
analytically isotopic to a point in $M\backslash E$.  But ${\cal
H}^{d-2}({\cal S}_{E,t,\tau,a,w_1^0})=0$. Therefore, it suffices to
shift slightly the parameter $s$ of $A_{t,\tau,a,\rho,w_1^0,s}$ in a
nearby parameter $s'$, which makes $A_{t,\tau,a,\rho,w_1^0,s}$ and
$A_{t,\tau,a,\rho,w_1^0,s'}$ isotopic to each other, so that ${\cal
D}_{t,\tau,a,w_1^0,s'}\cap E=\emptyset$.

The proof of Lemma 5.8 is complete.

\medskip
\noindent
{\it Step five: Holomorphic extension.}  Let $v_{00} \in T_0M_1$ with
$v_{00}\not\in T_0^cM$ and $|v_0-v_{00}|< \varepsilon$, let $C$ be a
$n$-dimensional proper linear cone in the $(2m+n)$-dimensional space
$T_0M$ and contained in ${\Gamma}_1$ such that $v_{00}\in C$, such
that the projection $T_0 C \to T_0 M / T_0^cM$ is surjective and such
that $\overline{C} \cap T_0^cM =\{0\}$.  Fix $p=p_1$ and set first
$p=0$. Note that $A_{t,\tau,a,\rho_1,0}(1)=0$. Let ${\cal P}_C$ denote
the set of parameters
$${\cal P}_C=\{(t,\tau,a)\in {\cal T} \times I_{\tau_0}\times {\cal
A}; \ \frac{d}{d\theta}A_{t,\tau,a}(1) \in C\},$$ which is a ${\cal
C}^1$-smooth $(n-1)$-dimensional submanifold of ${\cal T} \times
I_{\tau_0} \times {\cal A}$.  We choose a nearby piece of a manifold,
still denoted by ${\cal P}_C$, with same tangent space at 0 which is
${\cal C}^2$-smooth.

Let $K\subset M$ be a germ of a ${\cal C}^1$-smooth one-codimensional
submanifold of $M$ with $0\in K$ and $T_0 K \oplus \R
\frac{\partial}{\partial \theta} A(1) = T_0M$. Let ${\cal K}$ be small
neighborhood of $0$ in $K$. By $\Delta_1$, we denote a small
neighborhood of 1 in $\overline{\Delta}$ and
$\stackrel{\circ}{\Delta}_1 \subset \Delta$ its interior.

One observes that a consequence of the isotopy property 5.2 and of the
fact that the mapping
$${\cal P}_C\times \{\rho_1\} \times {\cal K} \times
\stackrel{\circ}{\Delta}_1 \ni (t,\tau,a,\rho_1,p,\zeta) \mapsto
A_{t,\tau,a,\rho_1,p}(\zeta) \in \C^{m+n}\backslash M$$ is a smooth
embedding is that ${\cal H}(\omega)$ extends holomorphically into the
open wedge set
$${\cal W}_{{\cal P}_C}=\{A_{t,\tau,a,\rho_1,p}(\zeta); \
(t,\tau,a)\in {\cal P}_C, p\in {\cal K}, \zeta \in
\stackrel{\circ}{\Delta}_1\}$$ minus the set
$$E_{{\cal P}_C}=\{A_{t,\tau,a,\rho_1,p}(\zeta); \ (t,\tau,a)\in {\cal
P}_C, p\in {\cal K}\cap \Phi, A_{t,\tau,a,\rho_1,p}(b\Delta)\cap E
\neq \emptyset, \zeta \in \stackrel{\circ}{\Delta}_1 \}.$$ Indeed,
since the mapping remains injective on ${\cal P}_C \times \{\rho_1\}
\times {\cal K}$ minus the set $\Phi_E\subset {\cal P}_C \times
\{\rho_1\} \times {\cal K}$ of $(t,\tau,a,\rho_1,\zeta)$ such that the
boundary disc meets $E$ at one or more points, we can set
unambiguously
$$F(z):= \frac{1}{2i\pi} \int_{b\Delta} \frac{f\circ
A_{t,\tau,a,p}(\eta)}{\eta - \zeta} d\eta$$ as a value at points
$z=A_{t,\tau,a,\rho_1,p}(\zeta)$ for an extension of $f|_{M\backslash
E}$, for such $(t,\tau,a,\rho_1,p)\in {\cal P}_C\times \{\rho_1\}
\times {\cal K} \backslash \Phi_E$.  Since $f$ extends holomorphically
to the interior of these discs, we get a continuous extension $F$ on
each $A_{t,\tau,a,\rho_1,p}(\Delta_1)$, $(t,\tau,a,\rho_1,p)\in {\cal
P}_C\times \{\rho_1\} \times {\cal K} \backslash \Phi_E$.  Thus, the
extension $F$ of $f|_{M\backslash E}$ also becomes continuous on
$$({\cal W}_{{\cal P}_C}\backslash {\cal E}_{\Phi_E}) \cup
(M\backslash E),$$ where ${\cal E}_{\Phi_E}$ is the proper closed
subset of ${\cal W}_{{\cal P}_C}$ defined by
$${\cal E}_{\Phi_E}=\{A_{t,\tau,a,\rho_1,p}(\zeta); \
(t,\tau,a,\rho_1,p)\in \Phi_E, \zeta \in \stackrel{\circ}{\Delta}_1
\}.$$ Since $f|_{M\backslash E}$ extends analytically to a
neighborhood of $A_{t, \tau, a, \rho_1, p}(\overline{\Delta})$, $F$ is
holomorphic in $({\cal W}_{{\cal P}_C}\backslash {\cal
E}_{\Phi_E})$. Indeed, fix a point $(\tilde{t}, \tilde{\tau},
\tilde{a}, \rho_1, \tilde{p}_0) \in {\cal P}_C \times \{\rho_1\}
\times ({\cal K} \backslash \Phi_E)$ and let $\tilde{\cal P} \times
\{\rho_1\} \times \tilde{\cal K}$ be a neighborhood of $(\tilde{t},
\tilde{\tau}, \tilde{a}, \rho_1, \tilde{p}_0)$ in ${\cal P}_C \times
\{\rho_1\} \times ({\cal K}\backslash \Phi_E)$ such that for each $(t,
\tau, a,\rho_1, p)\in \tilde{\cal P} \times \{\rho_1\} \times
\tilde{\cal K}$, $A_{t,\tau,a,\rho_1,p}(\overline{\Delta})$ is
contained in some neighborhood $\tilde{\omega}$ of $A_{\tilde{t},
\tilde{\tau}, \rho_1, \tilde{a}, \tilde{p}_0}(\overline{\Delta})$ in
$\C^{m+n}$ such that there exists a holomorphic function $\tilde{f}
\in {\cal H}(\tilde{\omega})$ with $\tilde{f}$ equal to $f$ near
$A_{\tilde{t}, \tilde{\tau}, \tilde{a}, \rho_1,
\tilde{p}_0}(b\Delta)$.  Let $\tilde{\zeta}\in
\stackrel{\circ}{\Delta}_1$ and $\tilde{z}=A_{\tilde{t}, \tilde{\tau},
\tilde{a}, \rho_1, \tilde{p}_0}(\tilde{\zeta})$. To check that the
previously defined function $F$ is holomorphic in a neighborhood of
$\tilde{z}$, we notice that for $z=A_{t, \tau, a, p}(\zeta)$, $(t, \tau,
a, \rho_1, p) \in \tilde{\cal P} \times \{\rho_1\}\times \tilde{\cal
K}$, $\zeta$ in some neighborhood $\tilde{\Delta}_1$ of
$\tilde{\zeta}$ in $\stackrel{\circ}{\Delta}_1$, $\tilde{f}(z)$ is
given by the Cauchy integral formula
$$\tilde{f}(z)=\frac{1}{2i\pi} \int_{b \Delta} \frac{\tilde{f} \circ
A_{t, \tau, a , \rho_1 ,p}(\eta)}{\eta-\zeta} d\eta= \frac{1}{2i\pi}
\int_{b \Delta} \frac{f\circ A_{t, \tau, a ,\rho_1
,p}(\eta)}{\eta-\zeta} d\eta= F(z).$$ As a consequence,
$\tilde{f}(z)=F(z)$ for $z$ in a small neighborhood of $\tilde{z}$ in
$\C^{m+n}$, since the mapping $(t, \tau, a, p,\zeta) \mapsto A_{t,
\tau, a, \rho_1, p}(\zeta)$ from $\tilde{\cal P}\times \tilde{\cal K}
\times \tilde{\Delta}_1$ to $\C^{m+n}$ has rank $2n$ at $(\tilde{t},
\tilde{\tau}, \tilde{a}, \tilde{p}_0, \tilde{\zeta})$.

This proves that $F$ is holomorphic into ${\cal W}_{{\cal P}_C}
\backslash {\cal E}_{\Phi_E}$.
 
By shrinking $\omega$ near $0$, which does not modify the possible
disc deformations, we can insure that $\omega \cap {\cal W}_{{\cal
P}_C}$ is connected, since $\overline{C}\cap T_0^cM=\{0\}$.  By Lemma
5.9 below the same is true for $\omega \cap {\cal W}_{{\cal P}_C}
\backslash {\cal E}_{\Phi_E}$.  Indeed, the fact that ${\cal
E}_{\Phi_E}$ is of zero $(2m+2n-1)$-dimensional Hausdorff measure
implies such connectedness.  Therefore $f\in {\cal H}(\omega)$ and
$F\in {\cal H} ({\cal W}_{{\cal P}_C} \backslash {\cal E}_{\Phi_E})$
stick together in a single holomorphic function in $\omega \cup ({\cal
W}_{{\cal P}_C} \backslash {\cal E}_{\Phi_E})$, since both are
continuous up to $M\backslash E$, which is a uniqueness set, and
coincide there.

\medskip
\noindent
{\sc 5.9. Lemma.} {\it The set ${\cal E}_{\Phi_E}$ is a union of
complex curves and ${\cal H}^{2m+2n-1}({\cal E}_{\Phi_E})=0$.}

\medskip
{\it Proof.} Namely, near the origin,
$${\cal E}_{\Phi_E}=\bigcup_{(t,\tau,a,\rho_1,p) \in \Phi_E}
A_{t,\tau,a,\rho_1,p}(\stackrel{\circ}{\Delta}_1).$$ Fix $t\in C$,
$\tau$ and $a$, consider the sets
$$M_{t,\tau,a}=\{A_{t,\tau,a,\rho_1,p}(\zeta); \ p\in {\cal K}, \zeta \in
\stackrel{\circ}{\Delta}_1\},$$ depending on $t$. Then ${\cal
W}_{{\cal P}_C}$ is foliated near $0$ by the manifolds $M_{t,\tau,a}$,
$\hbox{dim} \ M_{t,\tau,a} = 1 + \hbox{dim} \ M$, that contain $M$ in their
boundaries near $0$ and $M$ is foliated by the arcs
$A_{t,\tau,a,\rho_1,p}(\Delta_1\cap b\Delta)$. Notice that
$A_{t,\tau,a,\rho_1,p} (b\Delta)\cap E =A_{t,\tau,a\rho_1,p}(b\Delta \cap
\Delta_1)\cap E$ by the choice of $A=A_0$ with $A(1)\in M_1$, $A(b\Delta
\backslash \{1\}) \subset M_1^+\backslash M_1$, by $E \subset M_1^-$
and by continuity.  A direct application of Proposition 5.7 (iv)
entails that, for all fixed $(t,\tau,a)$, the set
$${\cal K}_{E,t,\tau,a}=\{p\in {\cal K}; \ A_{t,\tau,a,\rho_1,p}
(b\Delta \cap \Delta_1) \cap E \neq \emptyset\} \subset {\cal K},$$
contained in the $(d-1)$-dimensional manifold ${\cal K}$, satisfies
$${\cal H}^{d-2}({\cal K}_{E,t,\tau,a})=0.$$ Therefore, ${\cal
H}^{d}(M_{E,t,\tau,a})=0$ too in $M_{t,\tau,a}$ 
and since ${\cal W}_{{\cal P}_C}$ is
regularly foliated by the $M_{t,\tau,a}$,
$${\cal H}^{2m+2n-1}({\cal E}_{\Phi_E})={\cal H}^{2m+2n-1}
(\bigcup_{(t,\tau,a)\in {\cal P}_C} M_{E,t,\tau,a})=0.$$ The proof of
Lemma 5.9 is complete.

\smallskip
{\it Remark.} By Proposition 5.7 (iii), the set ${\cal K}\backslash
{\cal K}_{E,t,\tau,a}$ is connected. This provides another proof of
the isotopy property Lemma 5.8.

\smallskip
A closed set ${\cal E}$ contained in an open set ${\cal W}\subset
\C^{m+n}$ with ${\cal H}^{2m+2n-1}({\cal E})=0$ being not
automatically removable, we must study the structure of ${\cal
E}_{\Phi_E}$ with respect to the local complex structure of ${\cal
W}$.

\smallskip
\noindent
{\it The hypersurface case.}  Let us first give a proof for removing
${\cal E}_{\Phi_E}$ in the case $m=n=1$.  To begin the variations on
this theme, recall that ${\cal E}_{\Phi_E}$ would be removable if
${\cal H}^{2m+2n-2}({\cal E}_{\Phi_E})=0$, which
completes the proof of Theorem 5.1, (ii).

\medskip
\noindent
{\sc Theorem.} (see \cite{CH}). {\it Let ${\cal E}$ be a closed subset
of an open set ${\cal U} \subset \C^{m+n}$, $m+n\geq 1$, with ${\cal
H}^{2m+2n-2}({\cal E})=0$. Then, for every function $f\in {\cal
H}({\cal U} \backslash {\cal E})$, there exists a function $F\in {\cal
H}({\cal U})$ such that $F=f$ in ${\cal U} \backslash {\cal E}$.}

\medskip
{\it Proof.} We can, by localization, assume that ${\cal U}=P$ is the
unit polydisc $\Delta^{m+n}$ and that $0\in {\cal E}$. Let $G(k,m+n)$,
$0\leq k\leq m+n$, denote the grassmannian of $k$-dimensional complex
planes passing through the origin in $\C^{m+n}$.

\medskip
\noindent
{\sc Proposition.} (\cite{CH}). {\it Let ${\cal E}$ be a closed subset
in $\Delta^{m+n}$ such that ${\cal H}^{2k+1}({\cal E})=0$ for some
integer $k<m+n$. Then for almost every plane $L\in G(m+n-k,m+n)$,
${\cal H}^1({\cal E}\cap L)=0$.}

\medskip
Choose therefore a complex line $L$ through $0\in {\cal E}$ such that ${\cal
H}^1(L\cap \Delta^{m+n} \cap {\cal E})=0$ and an orthogonal complex
$(m+n-1)$-dimensional space $H$, $H\oplus L =\C^{m+n}$. By Proposition
(iv), for almost all $h\in H$, $\{h\} \times L\cap {\cal E}
=\emptyset$. Write $\Delta_L=L \cap \Delta^{m+n}$, choose a point
$p\in \Delta_L \cap {\cal E}$ ({\it e.g.} 0) and draw a small complex disc
$\Delta_L(p,r_0)$ of radius $r_0 >0$ and center $p$ with
$\overline{\Delta_L(p,r_0)} \subset \Delta_L$. For almost all $r<r_0$,
$b\Delta_L(p,r) \cap {\cal E} =\emptyset$, still thanks to Proposition
5.7 (iv). Furthermore, $\overline{\Delta_L(p,r)} + h \subset
\Delta^{m+n} \backslash {\cal E}$ for arbitrarily small $h\in
H$. Therefore, such a disc yields removability of $p$ for ${\cal
H}(\Delta^{m+n}\backslash {\cal E})$ along an obvious isotopy.

The proof is complete.

\medskip
\noindent
{\sc Proposition.} {\it Let $J$ and $K$ be closed subsets of $\Delta$
with ${\cal H}^1(J)={\cal H}^1(K)=0$ and set ${\cal E}=(J\times
\Delta) \cap (\Delta \times K)\subset \Delta^2$. Then ${\cal E}$ is
removable for ${\cal H}(\Delta^2 \backslash {\cal E})$.}

\medskip
{\it Proof.}
Fix $f\in {\cal H}(\Delta^2\backslash {\cal E})$.  For all $\zeta_2\in
\Delta \backslash K$, one has $\Delta \times \{\zeta_2\} \subset
\Delta^2 \backslash {\cal E}$. For $\zeta_2\in K$,
$$(\Delta \times \{\zeta_2\}) \cap {\cal E} =(\Delta \times
\{\zeta_2\}) \cap (J\times \Delta)= J \times \{\zeta_2\},$$ so ${\cal
H}^1(\Delta \times \{\zeta_2\})\cap E)=0$ and $f$ is holomorphic near
each point $(\zeta_1,\zeta_2) \in \Delta \times \{\zeta_2\} \backslash
{\cal E}$.  There exists $\lambda_2\in \C$ arbitrarily small with
$\Delta \times \{\zeta_2+\lambda_2\} \subset \Delta^2 \backslash {\cal
E}$.  Applying Proposition 5.7 (iv), we have ${\cal H}^0(rb\Delta
\times \{\zeta_2\} \cap {\cal E}) = 0$ for almost
all $0<r<1$. In other words, $rb\Delta \cap {\cal E}=\emptyset$. For
such $r$, the continuity principle along the family $(s,\zeta)\mapsto
(r\zeta,(1-s)\lambda_2+\zeta_2)$, $0\leq s \leq 1$, makes an extension
of $f$ in a neighborhood ${\cal V}(r\overline{\Delta} \times
\{\zeta_2\}) \subset \Delta^2$. Notice that ${\cal
V}(r\overline{\Delta}\times \{\zeta_2\}) \cap \Delta^2 \backslash
{\cal E}$ is connected, since $\Delta^2 \backslash (\Delta \times K)$
and $\Delta^2\backslash (J\times \Delta)$ are connected, by
Proposition 5.7 (iii), so uniqueness is guaranteed.

The proof is complete.

\smallskip
{\it Remark.} The argument above relies upon two facts: 1. Almost all
{\it complex} discs lay in $\Delta^2 \backslash {\cal E}$; 2. {\it
Every} disc touching ${\cal E}$ satisfies ${\cal H}^1(\Delta \cap
{\cal E})=0$.

\medskip
These preliminaries provide us with a proof of Theorem 5.1 in the case
that $M\subset \C^2$ is a ${\cal C}^{2,\alpha}$ hypersurface. Indeed,
the $a$-space in the maximal family of discs is empty and we can
choose small $(t_1,\tau_1)$ and $(t_2,\tau_2)$ so that $v_1$ and $v_2$
are linearly independent, $v_1=\partial/\partial \theta
A_{t_1,\tau_1,\rho_1,0}(1)$ and $v_2=\partial/\partial \theta
A_{t_2,\tau_2,\rho_1,0}(1)$ so that $v_1$ and $v_2$ point out in the
same side of $T_0^cM$ relative to $T_0M$ 
and $v_1,v_2 \not\in T_0^cM$. This is clearly
possible, by Lemma 5.6. Simply denote these two families by
$A_{1,p}(\zeta)$ and $A_{2,p}(\zeta)$. Then, since the normal bundle
to $M$ is of rank one and since $Jv_1, Jv_2$ point in the same side of
$M$ in $\C^2$, the two following wedges (same one-sided neighborhood)
with edge $M$ at $0$
$${\cal W}_1=\{A_{1,p}(\stackrel{\circ}{\Delta}_1); \ p\in {\cal K}\},
\ \ \ \ \ {\cal W}_2=\{A_{2,p}(\stackrel{\circ}{\Delta}_1); \ p\in
{\cal K}\}$$ contain a side ${\cal W}$ of boundary $M$ at $0$. As
above, one can construct two holomorphic extensions $F_1$ and $F_2$ of
$f$ to ${\cal W}_1 \backslash {\cal E}_1$ and ${\cal W}_2 \backslash
{\cal E}_2$ respectively, where ${\cal E}_1$ and ${\cal E}_2$ denote
the sets corresponding to the discs touching $E$:
$${\cal E}_1=\{A_{1,p}(\zeta); \ p\in {\cal K}, \zeta\in
\stackrel{\circ}{\Delta}_1, A_{1,p}(b\Delta \cap \Delta_1) \cap E \neq
\emptyset\}$$ and similarly for ${\cal E}_2$. By Proposition 5.7 (iv),
the one-dimensional Hausdorff measure of $\Phi_j=\{p\in {\cal K}; \
A_{j,p}(b\Delta)\cap E \neq \emptyset\}$, $j=1,2$, is zero. Indeed,
recall that $A(b\Delta)\cap E= A(b\Delta\cap \Delta_1) \cap E$ for
such discs with $p$ very close to 0 (in comparison with the size of
$\Delta_1$) and that $\{A_{j,p}(b\Delta \cap \Delta_1); \ p\in {\cal
K}\}$ foliate a neighborhood of $0$ in $M$.

Let $z\in {\cal W}\subset {\cal W}_1\cap {\cal W}_2$ such that $z\in
{\cal E}_1$ or $z\in {\cal E}_2$, say $z\in {\cal E}_1$,
$z=A_{1,p_{1,z}}(\zeta_{1,z})$ and $z=A_{2,p_{2,z}}(\zeta_{2,z})$.
Notice first that there are arbitrarily close to $p_{1,z}$ points
$p\in {\cal K}$ such that $A_{1,p}(b\Delta)\subset M\backslash E$.
Second, since $A_{1,p_{1,z}}$ and $A_{2,p_{2,z}}$ are transversal in
$\C^2$ at $z$ (by the choice of $(v_1,v_2)$), for all points $p$
varying in a neighborhood ${\cal V}(p_{1,z})\subset {\cal K}$, the
discs $A_{1,p_{1,z}} \cap A_{2,p} =\{z(p)\}$ intersect transversally
in a single point $z(p)$ such that $z(p_{2,z})=z$ and ${\cal V}(p_1)
\ni p\mapsto z(p)\in A_{2,p}( \stackrel{\circ}{\Delta}_1)$ is a local
${\cal C}^1$ diffeomorphism. But only for $p\in {\cal E}_2$ is the
disc $A_{2,p}$ not analytically isotopic to a point in $M\backslash E$
and ${\cal H}^3({\cal E}_2)=0$.  Thus, this shows that $F$ is already
holomorphic at each point of ${\cal V}(z) \cap
A_{1,p_{1,z}}(\stackrel{\circ}{\Delta}_1)$ outside a thin closed
subset $e_1 \subset A_{1,p_{1,z}}(\stackrel{\circ}{\Delta}_1)$ with
${\cal H}^1(e_1)=0$. Now, there exists a small circle contained in
this disc not meeting $e_1$ and the isotopy used in the proof above
can be applied once again in this analogous context to prove that $z$
is removable for $F\in {\cal H}({\cal W}\backslash ({\cal E}_1 \cap
{\cal E}_2))$.

\medskip {\it Remark.} A general proof in the hypersurface case ({\it
i.e.} for $m\geq 1$) follows along the same lines as above or by
reduction to $\C^2$ by slicing and using a separate analyticity
theorem like, for example, Shiffman's theorem below
(see Chirka and Stout, \cite{CS}, Section 4, for related reductions).
This completes the proof of Theorem 5.1 (i).

\medskip 
{\it End of proof of Theorem 5.1.} By construction, there exists $F\in
{\cal H}(\omega \cup ({\cal W}_{{\cal P}_C} \backslash {\cal
E}_{\Phi_E}))$ extending $f$. Let $A$ be a disc in the family
generating ${\cal W}_{{\cal P}_C}$ and assume that $A(b\Delta \cap
\Delta_1) \cap E \neq \emptyset$.  We shall remove
$A(\stackrel{\circ}{\Delta}_1)$.  Notice that, by construction,
$A(b\Delta\cap \Delta_1) \not\subset E$.

There always exist a point $p=A(\zeta_p)$, $\zeta_p \in
\stackrel{\circ}{\Delta}_1$, such that $p\in b\omega \cap {\cal
E}_{{\Phi}_E}$. Simply denote ${\cal W}={\cal W}_{{\cal P}_C}$ and
${\cal E}={\cal W}_{{\cal P}_C}\backslash (\omega \cup ({\cal
W}_{{\cal P}_C} \backslash {\cal E}_{\Phi_E}))$, $F\in {\cal H}({\cal
W}\backslash {\cal E})$.

By imitating the first step reduction, we can assume (see the
explanation below) that, after perharps changing $p$, there is a germ
of a one-codimensional manifold $M_1$ such that $p\in M_1$ and the
remaining part of ${\cal E}$ to be removed is contained in a
half-space $M_1^-$.

Indeed, a neighborhood ${\cal U}$ of $p$ in ${\cal W}$ is foliated by
discs of the family generating the wedge ${\cal W}$ ($={\cal W}_{{\cal
P}_C}$) and these analytic discs are integral real surfaces of a
subbundle, say $K$, of $T{\cal W}$, which they span. In this
neighborhood, one half of each disc lies in $\omega \cap {\cal U}$,
the other half is outside $\omega$ and all discs are transversal to
$b\omega$ (assuming, from the beginning, that $b\omega \backslash M$
is smooth after shrinking $\omega$).

Furthermore, since $A(b\Delta \cap \Delta_1) \not\subset E$ for each
disc of the family (hence $A(\stackrel{\circ}{\Delta}_1) \cap
\omega\neq \emptyset$), setting ${\cal E}={\cal W}_{{\cal P}_C}
\backslash (\omega \cup ({\cal W}_{{\cal P}_C} \backslash {\cal
E}_{\Phi_E}))$, one has that ${\cal W}$ and ${\cal W}\backslash {\cal
E}$ are $K$-minimal, {\it i.e.} both are a single $K$-orbit.

This is a key geometric fact ({\it cf.} \cite{JO3}, \cite{ME2},
\cite{PO1}, \cite{MP} and 5.1 here).

Then after introducing a similar set ${\cal A}$ as in Step one, and
after making use of the differential geometric lemma quoted in Step
one, we remove a point $p=A(\zeta_p) \in b\omega$ and the other ones
are removed similarly.

To complete the proof of Theorem 5.1 (iii), let us recall that if
$M$ is real analytic, it is known by works of Bloom-Graham or
Baouendi-Rothschild that there exists a real analytically parameterized
family of analytic discs attached to $M$ as $A_{t,\tau,a,\rho,p}$
filling a wedge at the base point so that the foliation
of ${\cal W}$ by pieces of 
$A(\stackrel{\circ}{\Delta}_1)$ 
is a real analytic foliation. Furthermore,
contrarily to the globally minimal case, it is superfluous
to deform $M$ step by step after removing points of $E$, since
1. $M$ is already minimal at every point and 2. The isotopy
lemma 5.8 holds without assuming that $E$ is contained in a half-space
$M_1^-$.

\medskip
\noindent
{\sc 5.10. Theorem.} {\it Let ${\cal U}\subset \C^{m+n}$ be a domain
$($connected$)$ equipped with a ${\cal C}^{\omega}$ foliation by complex analytic
curves and let ${\cal E}\subset {\cal U}$ be a closed subset which is
a union of leaves, with ${\cal H}^{2m+2n-1}({\cal E})=0$.  Then a
function $F\in {\cal H}({\cal U} \backslash {\cal E})$ extends
holomorphically in a neighborhood of a whole leaf $A$ whenever $F$
extends holomorphically through a single point of $A$.}

\medskip
{\it Remark.} When $m=n=1$, we provide a proof below that ${\cal E}$ is 
removable for a ${\cal C}^{2,\alpha}$ foliation in Theorem 5.13 below.

\medskip
{\it Proof.}  After the above reduction, the geometric assumption is:
there exists a neighborhood ${\cal U}$ of $p$ which is foliated by
complex analytic curves, ${\cal U}=\bigcup_{\theta\in D} A_{\theta}$,
$D\subset R^{2m+2n-2}$ small open set, there exists a closed set
${\cal G} \subset D$ with ${\cal H}^{2m+2n-3}({\cal G})=0$
(corresponding in the proof of Theorem 5.1 to the discs attached to
$M$ which meet $E$) and the set ${\cal E}$ to be removed (in ${\cal
U}\equiv {\cal W}$ near $p$) is a union of half discs, ${\cal
E}=\bigcup_{\theta\in {\cal G}} A_{\theta}^-$, where
$A_{\theta}^-=A_{\theta}\cap M_1^-$.

\medskip
\noindent
{\sc 5.11. Proposition.} {\it Let ${\cal U}\subset \C^{m+n}$ be a
small open set ${\cal C}^{\cal \omega}$-foliated by complex 
curves $A_{\theta}$,
${\cal U}=\bigcup_{\theta\in D} A_{\theta}$, $0\in {\cal U}$,
$D\subset R^{2m+2n-2}$ open, $0\in D$, let ${\cal G}\subset D$ be a
closed set with ${\cal H}^{2m+2n-3}({\cal G})=0$, let $M_1$ be a
${\cal C}^1$ hypersurface through $0$ in ${\cal U}$ with $T_0M_1 +
T_0A_0=T_0\C^{m+n}$ and set ${\cal E}= (\bigcup_{\theta\in {\cal G}}
A_{\theta})\cap M_1^-$.  Then there exists a neighborhood ${\cal V}$
of $0$ such that for every function $f\in {\cal H}({\cal U}\backslash
{\cal E})$, there exists a function $f\in {\cal H}({\cal V})$ with
$F=f$ in ${\cal V}\backslash {\cal E}$.}

\medskip
{\it Proof.} It is {\it not} true that the foliation of ${\cal W}$ by
discs of the foliation is a complex analytic foliation, {\it i.e.}
locally equivalent to $\C\times \C^{m+n-1}$ after a
biholomorphism. Neither in 5.10, 5.11.

Nevertheless, let us first investigate geometrically this case.

\smallskip
\noindent
{\it The case of a holomorphic foliation.}  Here, the geometric
situation is that there exists a smooth hypersurface $M_1$ through $0$
and a closed set ${\cal G} \subset \Delta^{m+n-1}$ with ${\cal
H}^{2m+2n-3}({\cal G})=0$ such that ${\cal E}=M_1^-\cap (\Delta \times
{\cal G})$ near $0$. Indeed, one just straightens the holomorphic
foliation.

Notice that by Proposition 5.7 (iv), for almost all two-dimensional
affine complex planes $L\equiv \C\times \C$, $L\cap (\C \times {\cal
G})=\C\times {\cal G}_L\subset \C\times \C$ for a closed set ${\cal
G}_L \subset \C$ with ${\cal H}^1({\cal G}_L)=0$. Hence we are in the
following situation (a particular case of Proposition 5.11).

\medskip
\noindent
{\sc 5.12. Lemma.} {\it Let ${\cal U}$ be a connected open set in
$\C^2_{w,z}$, $0\in {\cal U}$, let $M_1\subset {\cal U}$ be a closed
hypersurface, $0\in M_1$, $T_0M_1 \oplus \R_u =T_0\C^2$, $u=\hbox{Re}
\ w$, and let ${\cal E}= (\C_w \times E) \cap {\cal U} \cap M_1^-$ be
closed, where $E\subset \C_z$ is closed and ${\cal
H}_{loc}^1(E)=0$. Then there exist a neighborhood ${\cal V}$ of $0$
such that for every function $f\in {\cal H}({\cal
U}\backslash {\cal E})$, there exists a
function $F\in {\cal H}({\cal V})$ with $F=f$ in ${\cal
V}\backslash {\cal E}$.}

\medskip
{\it Proof.} Notice that ${\cal H}^3({\cal E})=0$.  Define
$B_w(\zeta)=(w,r\zeta)$, $|w|$ small and $r>0$. By the fact that
${\cal H}^1(E)=0$, then for almost all $r>0$, the boundary of the disc
$\zeta \mapsto r\zeta$ does not meet $E$. Hence also $B_w(b\Delta)
\cap {\cal E}=\emptyset$ for such $r>0$, {\it because} ${\cal
E}\subset \C_w \times E$. Fix such $r$. Then all $B_w$ for different
$w$ are analytically isotopic to each other {\it in} ${\cal
U}\backslash {\cal E}$, $B_0(0)=0$ and, moreover, $B_u$ is
analytically isotopic to the point $(u,0)$ in ${\cal U}\backslash
{\cal E}$ if $u>0$. Therefore $0$ is removable.

The proof of Lemma 5.12 is complete.

\smallskip
For general $m+n\geq 2$, the above constructed isotopies
lay inside a fixed complex plane $L$, so that the continuity
principle in $\C^{m+n}$ applies, giving holomorphic extension at 0.

\smallskip
{\it Remark.} The ${\cal C}^{2,\alpha}$ (even ${\cal C}^1$) foliated version of Theorem 5.10
admits a proof in $\C^2$ that we give below.

\medskip
\noindent {\sc 5.13. Theorem.}
{\it Let ${\cal U}\subset\C^2$ be a domain equipped with a
$C^1$ foliation $\cal F$ by complex analytic curves. 
Further let ${\cal
E}$ be a closed union of leaves with ${\cal H}^3({\cal E})=0$. If
a function $f\in{\cal H}(D-{\cal E})$ admits a holomorphic
extension into the neighborhood of some point $p\in {\cal E}$, then it
extends analytically into a neighborhood of the leaf $L$ containing
$p$.}

\medskip
{\it Proof.}  Let $L'$ be the set of all points $z\in L$ such that $f$
extends holomorphically into a neighborhood of $z$. As all such
extensions obviously fit together, it is enough to show $L'=L$.

Let us assume $L'\not= L$. Clearly $L'$ is an open subset of $L$ with
respect to the leaf-topology. Fix a point $q$ on the relative boundary
of $L'$. In a neighborhood $U$ of $q$ we choose a holomorphic function
$z$ with $z(q)=0,dz(q)\not=0$, such that the curve $\{z=0\}$ intersects
the leaf $L$ transversely in $q$ . Hence 
${\cal H}^1({\cal E}\cap\{z=0\})=0$ (after contraction of $U$), 
and we may draw a simple closed curve $\gamma_0\subset\{z=0\}$ 
surrounding $q$ and avoiding ${\cal E}\cap\{z=0\}$.

If we choose $\gamma_0$ in a small neighborhood of $q$, we get a family
of contours $\gamma_{\zeta}\subset\{z=\zeta\}-{\cal E},|\zeta|\leq\epsilon$ 
by flowing ({\it i.e.} moving) $\gamma_0$ along the foliation. 
By hypothesis, for an open set of parameters $\zeta$,
the restriction of $f$ to $\{z=\zeta\}$ is holomorphic near the closure
of the domain $G_\zeta$ surrounded by $\gamma_\zeta$. 
As $f$ is holomorphic near $\bigcup_{|\zeta|\leq\epsilon}\gamma_\zeta$,
it extends by the continuity principle to a neighborhood of
$\bigcup_{|\zeta|\leq\epsilon}\overline{G}_\zeta$, in contradiction to
the choice of $q$.

The proof of Theorem 5.13 is complete.

\smallskip
{\it Remark.} However, the reduction to $\C^2$ is impossible
for a general ${\cal C}^{2,\alpha}$ foliation, since there
need not be families of complex surfaces foliated by complex
curves of the foliation. Therefore, 5.13 in $\C^2$ does not provide
5.13 for any $\C^{m+n}$. 

\medskip
\noindent
{\it The case of a real analytic foliation.} 
Here, we shall use a separate analyticity theorem due to Shiffman.
A subset $Q$ of a polydisc $\Delta^{m+n-1}$ is said to be a {\it full
subset} of $\Delta^{m+n-1}$ if $D\cap Q$ is a set of full measure in
$D$ for almost every coordinate disc $D\subset \Delta^{m+n-1}$.

\medskip
\noindent
{\sc Theorem.} ({\sc Shiffman}, \cite{SH}). {\it Let
$\Delta^{m+n-1}\subset \subset \C^{m+n-1}$ be a polydisc and let $Q
\subset \Delta^{m+n-1}$ be a full subset of $\Delta^{m+n-1}$. Then a
function $F: Q \to \C$ has a holomorphic extension to $\Delta^{m+n-1}$
if and only if, for almost every coordinate disc $D \subset
\Delta^{m+n-1}$, $F|_{D\cap Q}$ extends holomorphically to $D$.}

\medskip
First, $p=0$ in coordinates $(w,z)
\in \C\times \C^{m+n-1}$, $w=u+iv$, with $\C\times \{0\}=T_0A_0$ and
$T_0M_1=\{u=0\}$, $T_0M_1^-=\{u\leq 0\}$.  There are, in the whole the
grassmannian of affine complex lines passing near $0$ in ${\cal U}$,
the lines $a+h(L_0)$, $a\in \{0\}\times \C^{m+n-1}$ close to $0$, with
$h\in GL(m+n,\C)$ close to $Id$, where $L_0=\C_w\times 0$ which are
cut by $M_1$ in two pieces. Draw a small enough analytic disc
$B(\zeta)=(c(\zeta-1)-b,0)$, $\zeta\in \Delta$, inside $L_0$, with
$c>0$ small and fixed throughout, with $b>0$, $b<<c$ and define
$B_{a,h}(\zeta):=a+h\circ B(\zeta)$ so that $B_{a,h}(\Delta) \subset
a+h(L_0)$.  Notice that for small $|a| <<c$, $||h-Id||<<c$, then
$B_{a,h}(b\Delta)\cap M_1^-\subset B_{a,h}(b\Delta \cap \Delta_1)$,
for fixed $\Delta_1= \{|\zeta-1|\leq c_1\}\cap \overline{\Delta}$,
with $c_1=c/5$ say.  (In other words, the boundaries
$B_{a,h}(b\Delta)$ can meet ${\cal E}$ only along a fixed part of
their boundaries.) And notice that for fixed $h,c$ then $\bigcup_{a}
B_{a,h}(\Delta)$ makes a holomorphic foliation by complex discs with
the origin point in its interior.  Varying $h$, to apply Shiffman's
theorem, it suffices to show that for almost all $a\in \C^{m+n-1}$
close to $0$, a fixed function $F\in {\cal H}({\cal U}\backslash {\cal
E})$ extends holomorphically to $B_{a,h}(\Delta)$.

Clearly, if $c$ is small, $B_{a,h}(b\Delta)\cap M_1^-\subset
B_{a,h}(b\Delta \cap \{|\zeta -1|\leq 5b\})$ for all $a,h$, so that
$B_{a,h}(b\Delta\cap \Delta_1)$ is much longer than its intersection
with $M_1^-$.

Let $\Sigma_1 =\bigcup_a B_{a,h}(b\Delta \cap \Delta_1)$.  The
hypersurface $\Sigma_1$ is transversal to the real foliation, hence
${\cal H}^{2m+2n-2}(\Sigma_1 \cap {\cal E})=0$. Hence by Proposition
5.7 (iv), for almost all $a$, $B_{a,h}(b\Delta \cap \Delta_1) \cap
{\cal E}=\emptyset$ (by the property ${\cal H}^1(B_{a,h}(b\Delta\cap
\Delta_1))=0$). In other words, for almost all $a$,
$B_{a,h}(b\Delta)\cap {\cal E}=\emptyset$. To complete the proof of
Theorem 5.10, the remaining point is to isotope the boundaries
$B_{a,h}(b\Delta)$ analytically to a point inside ${\cal U}\backslash
{\cal E}$.

By the foliation assumption, there exist some {\it real analytic} coordinates
$(u,v,x)$ so that ${\cal U}=\R_{u,v}^2\times \R^{2m+2n-2}_x$ (here,
$(u,v)$ are real coordinates in general distinct from $(\hbox{Re} \ w,
\hbox{Im} \ w)$) such that the sets $\R_{u,v}^2\times \{ct.\}$
correspond to discs of the foliation and ${\cal E}=M_1^-\times
(\R^2_{u,v} \times {\cal G})$. We can assume that the Jacobian matrix
of the change from the holomorphic coordinates $(w,z)$ to $(u,v,x)$ is
the identity at the origin (so that $u$ and $v$ are close to
$\hbox{Re} \ w$ and $\hbox{Im} \ w$).

Denote $\gamma_B:= B_{a,h}(b\Delta \cap \Delta_1)$ and set
$$\Sigma_{a,h}=\{p+u\in {\cal U}; \ p\in \gamma_B, -5b<u<5b\}$$ in
these real analytic coordinates and come back to the original holomorphic
coordinates.

Then $\Sigma_{a,h}$ is a piece of a ${\cal C}^{\omega}$-smooth surface near
$0$ with the property that $\Sigma_{a,h} \cap {\cal E}=\emptyset$.
However, $\Sigma_{a,h}$ is not immediately seen to be a union of boundaries of
analytic discs as in the case of a holomorphic foliation.

Therefore, we shall apply a complexification argument.

From the beginning, we can assume that $b<<c$, $|a|<<c$, $|a|\approx
b$ and $||h-Id||<<|a|$ also, since the set of such $h$ still contain
an open set in $\hbox{GL}(m+n,\C)$, which is sufficient to apply
separate analyticity locally. Therefore, for all $p\in \Sigma_{a,h}$,
the tangent space $T_p\Sigma_{a,h}$ (considered as a linear subspace
of $T_0\R^{2m+2n}$ is close to $T_0 L_0=\C_w \times \{0\}$. Therefore,
$\Sigma_{a,h}$ is a graph over a domain $D_{a,h}$ contained in the
$w$-space. Since the transformation from $(w,z)$ to $(u,v,x)$ is close
to the identity in the ${\cal C}^1$ norm and since $||h-Id||$ and
$|a|$ are very small, this domain $D_{a,h}\subset \C_w$ is
approximately the domain $$D:=\{c(\zeta-1)-b+u_0 \in \C_w; \ \zeta\in
b\Delta, |\zeta-1|<c/5, -5b<u_0< 5b\}.$$ Anyway, by taking
$\Sigma_{a,h}$ a little bit larger, $D_{a,h}$ will contain $D$ for all
small $|a|, ||h-Id||$.

Hence there exist ${\cal C}^{\omega}$ functions $s_j: D\to \C$, $1\leq j\leq
m+n-1$, such that $\Sigma_{a,h}$ is given by the equations
$z_1=s_1(u,v),...,z_{m+n-1}=s_{m+n-1}(u,v)$ as a graph over $D$.

Notice that the domain $D$ is foliated by the real analytic arcs
$\gamma_{u_0}: b\Delta \cap \Delta_1 \ni \zeta \mapsto
c(\zeta-1)-b+u_0\in D$.

We write them as $[-d,d]\ni \theta \mapsto (u_{u_0}(\theta),
v_{u_0}(\theta))$, with $\zeta=e^{i\theta}$ and $c/5=|e^{id}-1|$.
Notice that by the disposition of $M_1$, if $\pi_w: \C^{m+n} \to \C_w$
denotes the projection, then $\pi_w(M_1^-) \cap
\gamma_{u_0}([-d,d])=\emptyset$ for $2b\leq u_0\leq 5b$. Hence
${\Sigma}_{a,h}$ is foliated by the analytic arcs $\Gamma_{u_0}: \theta
\mapsto (\gamma_{u_0}(\theta),
s_1(u_{u_0}(\theta),v_{u_0}(\theta)),...,
s_{m+n-1}(u_{u_0}(\theta),v_{u_0}(\theta)))$ and for $u_0 \geq 2$,
these arcs are far from $M_1^-$.

Now, we complexify $\theta$ in a complex variable $\Theta\in
[-d,d]+i[-5b,5b]=: T$ such that $\hbox{Re} \ \Theta=\theta$. Since
$u_{u_0}$, $v_{v_0}$ and the $s_j$ are analytic, such a complexification
exists and $\gamma_{u_0}(\Theta)$ makes a biholomorphism from this
strip $T$ to a strip neighborhood of $\gamma_{u_0}([-d,d])$ in $D$
which contains, say, $\{\gamma_{u_0}(\theta)+u; \ \theta \in [-d,d],
-3b \leq u\leq 3b\}=:D_{u_0}$. The domains of definition
of the complexification are uniform since the piece $\Sigma_{a,h}$
comes from the real analytic foliation, which can be supposed
to be given by converging series in a fixed neighborhood of 0. 
By the implicit function theorem over
$D_{u_0}$, we can replace $\Gamma_{u_0}(\Theta)$ by a parameterizing
variable $w=\Gamma_{u_0}(\Theta)$ 
so that $\Gamma_{u_0}( \stackrel{\circ}{T})$ will be a
complex manifold of dimension one given by a graph $D_{u_0}\ni w
\mapsto (w,\varphi_1(w),...,\varphi_{m+n-1}(w))\in \C^{m+n}$ and which
we will denote by $E_{u_0}$. Notice that
$E_{u_0}\subset\subset {\cal U}$. Therefore, for all $u_0$ with
$0\leq u_0<5b$, $E_{u_0}$ intersects $M_1$ transversally and
$\Gamma_{u_0} ([-d,d])\cap M_1^-\subset \Gamma_{u_0}([-5b,5b])$.

Now, {\it inside} the complex curve $E_{u_0}$, we can close-up the analytic arc
$\Gamma_{u_0}([-d,d])$ outside $M_1^-$, making the boundary of an
analytic disc $B_{u_0}$ parameterized by $u_0$, a part of its boundary
being given by $\gamma_{u_0}([-d,d])$ and the other part living in
$M_1^+\backslash M_1$.

Notice that, by construction, all the $B_{u_0}$ for $0\leq u_0\leq 2b$
are analytically isotopic to each other and that all their boundaries
are contained in ${\cal U}\backslash {\cal E}$. But $\Gamma_0(D)$ is
a piece of the disc $B_{a,h}(\Delta)$ (recall that $\gamma_B$ is
contained in ${\Sigma}_{a,h}$), so that $B_0$ is analytically isotopic to
$B_{a,h}$ in ${\cal U}\backslash {\cal E}$ (isotope their boundaries
inside $a+h(L_0) \backslash M_1^-$). And $B_{2b}$ is analytically
isotopic to a point in $M_1^+\backslash M_1$.  This yields the desired
isotopy.

The proofs of Theorems 5.10 and 5.1 are complete.

\end{document}